\documentclass[11pt,reqno]{amsart}
\usepackage[a4paper,margin=3cm]{geometry}
\usepackage{amsmath,amsfonts,amssymb,bm,graphicx,amsthm,color,yfonts,cite, latexsym}
\usepackage[hidelinks]{hyperref}
\usepackage{cite}
\usepackage{booktabs}
\usepackage{array}
\usepackage{tabularx}
\usepackage{float}
\usepackage{hyperref}
\usepackage{amsfonts}
\usepackage{empheq}
\usepackage{enumerate,mathtools,bbm,bm,tcolorbox,mathrsfs}
\mathtoolsset{showonlyrefs}
\usepackage{ragged2e}

\usepackage{tikz}
\usepackage{tikz-cd}
\usetikzlibrary{cd}
\usetikzlibrary{arrows.meta}
\usetikzlibrary{matrix, arrows, positioning}
\usetikzlibrary{shadows}
\usetikzlibrary{arrows.meta}
\colorlet{diagramcolor}{black}
\usetikzlibrary{calc}
\usepackage{caption}
\allowdisplaybreaks

\newcommand{\R}{\mathbb{R}}

\newcommand{\N}{\mathbb{N}}
\newcommand{\Z}{\mathbb{Z}}
\newcommand{\interior}{\operatorname{int}}
\newcommand{\supp}{\operatorname{supp}}
\newcommand{\Tr}{\operatorname{Tr}}
\newcommand{\one}{\mathbbm{1}}
\newcommand{\rank}{\operatorname{rank}}

\newtheorem{theorem}{Theorem}[section]
\newtheorem{proposition}[theorem]{Proposition}
\newtheorem{lemma}[theorem]{Lemma}

\newtheorem{question}[theorem]{Question}
\theoremstyle{remark}
\newtheorem{remark}[theorem]{Remark}
\numberwithin{equation}{section}

\makeatletter
\def\l@section{\@tocline{1}{0pt}{0pt}{}{}}%
\def\l@subsection{\@tocline{2}{0pt}{2.5em}{}{}}%
\def\l@subsubsection{\@tocline{3}{0pt}{3.0em}{}{}}%
\makeatother

\title[Sharp orthonormal Strichartz estimates]{Sharp orthonormal Strichartz estimates for the Schr\"odinger equation in all remaining cases}

\author[Y. Gong]{Yinghao Gong}
\address{Nagoya University}
\email{gong.yinghao.m7@s.mail.nagoya-u.ac.jp}

\author[S. Hadama]{Sonae Hadama}
\address{The University of Osaka}
\email{hadama.sonae.sci@osaka-u.ac.jp}

\author[S. Toyoshima]{Shunya Toyoshima}
\address{Saitama University}
\email{s.toyoshima.845@ms.saitama-u.ac.jp}

\keywords{Strichartz estimates, Orthonormal Strichartz estimates}
\subjclass[2020]{Primary 35Q41; Secondary 35B45, 42B37.}
\thanks{The authors are deeply grateful to Professor Rupert L. Frank for his valuable advice on the first version of this paper posted on arXiv. S.H. was supported by JSPS KAKENHI Grant Number JP26KJ0203.}

\date{}

\begin{document}
	
\begin{abstract}
	We study orthonormal Strichartz estimates for the free Schr\"odinger equation in homogeneous Sobolev spaces $\dot H^s$. We resolve the remaining open cases concerning the optimal summability of the coefficients. Together with earlier results and counterexamples, our theorems complete the classification for nonnegative Sobolev regularity $0\le s<d/2$.
\end{abstract}
	
	\maketitle
	\tableofcontents
	
	\section{Introduction}
	Let $U(t) = e^{it\Delta_x}$. 
	The classical \textit{Strichartz estimate} for the free Schr\"odinger equation with initial data $f\in\dot{H}^s_x(\R^d)$ states that 
	\begin{equation}\label{eq:single Strichartz}
		\| |U(t)f|^2 \|_{L_t^p(\R;L_x^q(\R^d))}\lesssim\|f\|_{\dot{H}_x^s(\R^d)}^2
	\end{equation}
	holds for all spatial dimensions $d\geq1$ if the exponents $p,q\in [1,\infty]$ and $s\in [0,d/2)$ satisfy
	\begin{equation}\label{eq:exponents}
		\frac2p+\frac dq = d-2s, \qquad (p,q)\ne (1,\infty).
	\end{equation}
	This type of estimate was first obtained by Strichartz via Fourier restriction \cite{Strichartz 1977}.
	Ginibre and Velo \cite{Ginibre Velo 1985, Ginibre Velo 1992} extended it to a wider range of exponents.
	Finally, Keel and Tao \cite{Keel Tao 1998} proved the endpoint case $p=1$ when $d\ge 3$ and consequently obtained \eqref{eq:single Strichartz} for all admissible exponents $p$ and $q$ when $s=0$. 
	By the Sobolev embedding and \eqref{eq:single Strichartz} with $s=0$, we obtain \eqref{eq:single Strichartz} for all $0\le s < d/2$ when $q<\infty$.
	Later, Guo, Li, Nakanishi, and Yan \cite[Proposition 2.1]{Guo et al 2018} proved \eqref{eq:single Strichartz} when $q=\infty$ and $1<p<\infty$.  
	
	One of the most natural questions is whether we can generalize classical Strichartz estimates to orthonormal initial data.
	If one naively uses \eqref{eq:single Strichartz} and the triangle inequality, one can get only
	\begin{equation}\label{eq:naive}
		\bigg\|\sum_j\nu_j|U(t)f_j|^2\bigg\|_{L_t^p(\R;L_x^q(\R^d))}\lesssim\|\nu\|_{\ell^1}
	\end{equation}
	for all orthonormal systems $(f_j)_j$ in $\dot{H}^s_x(\R^d)$ and all sequences $\nu=(\nu_j)_j$.
	Note that we can prove this without using the orthogonality of $(f_j)_j$.
	It has been revealed recently that, for suitable choices of $p,q,s$, we can improve \eqref{eq:naive} to    
	\begin{equation}\label{eq:OS}
		\bigg\|\sum_j \nu_j |U(t)f_j|^2\bigg\|_{L^p_t(\R;L^q_x(\R^d))} \lesssim \|\nu\|_{\ell^\alpha}
	\end{equation}
	for some $\alpha>1$ by using the orthogonality of $(f_j)_j$.
		
		The study of \textit{orthonormal Strichartz estimates} was initiated by	Frank, Lewin, Lieb, and Seiringer \cite{Frank et al 2014}.
		For $s=0$, they established such estimates with the optimal summability exponent $\alpha=2q/(q+1)$ for a range of exponents $p,q$.
		Frank and Sabin \cite{Frank Sabin 2017} extended the range of $p,q$ for which this optimal $\alpha=2q/(q+1)$ is attained.
		 Bez, Hong, Lee, Nakamura, and Sawano \cite{Bez et al 2019} dealt with the case $0<s<d/2$. Some boundary and endpoint cases were subsequently studied in \cite{Bez et al 2020,Bez et al 2024}. For more recent developments, we refer the reader to \cite{Hadama Toyoshima 2026+, Gong 2026-, Gong 2026- 2, Feng et al 2026-}.
		See Sections \ref{subsec:main} and \ref{subsec:history} for a more detailed account of previous results.
		
		We can regard previous studies that were mentioned above as efforts to answer the following question.

	\begin{question}\label{question}
		Let $p,q\in [1,\infty]$ and $s\in[0,d/2)$ be such that \eqref{eq:exponents} holds.
		There are two possible cases: there exists $\alpha_*\in[1,\infty]$ such that
		\begin{enumerate}[$(i)$]
			\item The estimate \eqref{eq:OS} holds if $\alpha \le \alpha_*$ but fails if $\alpha > \alpha_*$;
			\item The estimate \eqref{eq:OS} holds if $\alpha < \alpha_*$ but fails if $\alpha \ge \alpha_*$.
		\end{enumerate}
		The first question is to determine $\alpha_*$.
		The second question is which of $(i)$ and $(ii)$ occurs.
		In case $(ii)$, we further ask whether the estimate
		\begin{equation}\label{eq:OS restricted}
			\bigg\|\sum_j \nu_j |U(t)f_j|^2\bigg\|_{L^p_t(\R;L^q_x(\R^d))} \lesssim \|\nu\|_{\ell^{\alpha_*,\beta}}
		\end{equation}
		holds for some $1\le \beta < \alpha_*$, where $\ell^{\alpha,\beta}$ is the Lorentz space (see Section \ref{subsec:notation}).
	\end{question}
	
	Despite significant progress on this topic in recent years, several cases of Question \ref{question} remained open.
	In this paper, we resolve \textit{all} remaining cases.

	\subsection{Main results}\label{subsec:main}
	When $d=1$, the only remaining case is $(p,q)=(2,\infty)$, namely, the point $A$ in Figure \ref{fig:1}.
	Bez, Lee, and Nakamura \cite{Bez et al 2020} proved that a weaker version
	\begin{equation}\label{eq:weaker version}
	\bigg\|\sum_j \nu_j |U(t)f_j|^2\bigg\|_{L^{2,\infty}_t (\R;L^\infty_x(\R))}
	\lesssim \|\nu\|_{\ell^\alpha}
	\end{equation}
	holds for all $1\le \alpha <2$. 
	Moreover, they proved that we can improve the $L^{2,\infty}_t$-norm in \eqref{eq:weaker version} to the $L^2_t$-norm when $1\le \alpha \le 4/3$.
	On the other hand, via a Kakeya set, it was proved that the critical summability estimate \eqref{eq:OS restricted} fails for $\alpha_* =2$ and for all $\beta \ge 1$ in \cite{Bez et al 2019}.
	Therefore, it remained open whether \eqref{eq:OS} holds when $4/3<\alpha <2$.
	Our first main result claims that indeed it holds.
	\begin{theorem}[At \text{$A$} when \text{$d=1$}]\label{thm:d=1}
		Let $d=1$, $p=2$, and $q=\infty$.
		Then, for any $4/3< \alpha<2$, the estimate \eqref{eq:OS} holds for every $\nu\in \ell^\alpha$ and for every orthonormal system $(f_j)_j$ in $L^2_x(\R)$.
	\end{theorem}
	
\begin{figure}[!b]
	\centering
			\begin{minipage}[b]{4cm}
		\centering
		\begin{tikzpicture}
			\draw[black,thin,densely dotted,-stealth] (0,0)node[below left]{$O$}--(3,0)node[below]{$1/q$};
			\draw[black,thin,densely dotted] (5/2,0)--(5/2,5/2)--(0,5/2);
			\draw[black,thin,densely dotted] (0,0)--(5/2,5/2);
			\draw[black,thin,densely dotted,-stealth] (0,5/4)--(0,3)node[left]{$1/p$};
			\draw[black,thin,densely dotted] (0,5/4)node[above right]{$A=C=D$}node[left]{$1/2$}--(5/2,0)node[above right]{$B$}node[below]{$1$};
			\draw[black,thin,densely dotted] (0,0)--(0,5/4);
			\filldraw[fill=white] (0,0) circle[radius=0.06];
			\filldraw[fill=black] (0,5/4) circle[radius=0.06];
			\node[right] at (4/7,4/8) {$F$};
		\end{tikzpicture}
		\caption*{$d=1$}
	\end{minipage}
	\hspace{+6mm}
	\begin{minipage}[b]{4cm}
		\centering
		\begin{tikzpicture}[scale=1.0]
			\draw[black,thin,densely dotted,-stealth] (0,0)node[below left]{$O$}--(3,0)node[below]{$1/q$};
			\draw[black,thin,densely dotted,-stealth] (0,5/2)--(0,3)node[left]{$1/p$};
			\draw[black,thin,densely dotted] (5/2,0)--(5/2,5/2)--(0,5/2);
			\draw[black,thin,densely dotted] (5/2,0)node[above right]{$B$}node[below]{$1$}--(5/6,5/3);
			\draw[black,thin,densely dotted] (0,0)--(5/2,5/2);
			\draw[black,thin,densely dotted] (0,5/4)node[left]{$1/2$}--(5/2,5/2);
			\draw[black,very thick,solid] (0,0)--(0,5/2);
			\draw[black,very thick,solid] (5/6,5/3)--(0,5/2);
			\draw[black,very thick,solid] (0,0)--(5/6,5/3);
			\filldraw[fill=white] (0,0) circle[radius=0.06];
			\filldraw[fill=white] (0,5/2) circle[radius=0.06]node[above right]{$C=D$}node[left]{$1$};
			\node[above right] at (5/6,5/3) {$A$};
			\node[right] at (5/4,5/4) {$F$};
		\end{tikzpicture}
		\caption*{$d=2$}
	\end{minipage}
	\hspace{+6mm}
	\begin{minipage}[b]{4cm}
		\centering
		\begin{tikzpicture}[scale=1.0]
			\draw[black,thin,densely dotted,-stealth] (0,0)node[below left]{$O$}--(3,0)node[below]{$1/q$};
			\draw[black,thin,densely dotted,-stealth] (0,5/2)--(0,3)node[left]{$1/p$};
			\draw[black,thin,densely dotted] (5/2,0)--(5/2,5/2)--(0,5/2);
			\draw[black,thin,densely dotted] (0,0)--(5/2,5/2);
			\draw[black,thin,densely dotted] (0,5/4)node[left]{$1/2$}--(5/2,5/2);
			\draw[black,thin,densely dotted] (3/2,5/2)node[above]{$D$}--(5/2,0)node[above right]{$B$}node[below]{$1$};
			\draw[black,very thick,solid] (0,0)--(5/3,25/12);
			\draw[black,very thick,solid] (0,0)--(0,5/2);
			\filldraw[fill=white] (0,0) circle[radius=0.06];
			\filldraw[fill=white] (0,5/2) circle[radius=0.06]node[above right]{$C$}node[left]{$1$};
			\node[above right] at (5/3,25/12) {$A$};
			\node[right] at (25/14,25/14) {$F$};
		\end{tikzpicture}
		\caption*{$d\geq3$}
	\end{minipage}
	\caption{The thick black segments and the filled point indicate the cases resolved in this paper. The endpoints of the segments are excluded. Note that $A=(\frac{d-1}{d+1}, \frac{d}{d+1})$.}
	\label{fig:1}
\end{figure}
	When $d=2$, the remaining cases are the segments $(O,A)$, $(A,D)$, and $(O,C)$.
	The segments $(O,A)$ and $(O,C)$ also arise in higher dimensional cases, but $(A,D)$ remains open only when $d=2$.
	The second and third authors \cite{Hadama Toyoshima 2026+} proved the endpoint estimates at $A$ when $d\ge 2$.
	When $d\ge 3$, the point $D$ is Keel--Tao's endpoint. Hence, real interpolation between $A$ and $D$ gives optimal estimates on $(A,D)$. However, the endpoint Strichartz estimate fails when $d=2$ (see \cite{Montgomery-Smith 1998}), and the optimal estimates on $(A,D)$ remained open. 
	Our second main result solves this problem.
	\begin{theorem}[On \text{$(A,D)$ when $d=2$}]\label{thm:d=2}
		Let $d=2$, $1<p<3/2$ and $q'=p$.
		Then, for every orthonormal system $(f_j)_j$ in $L^2_x(\R^2)$ and for every sequence $\nu=(\nu_j)_j\in\ell^{p}$, the estimate \eqref{eq:OS} holds with $\alpha=p$.
	\end{theorem}

	When $d\ge 2$, there are now two remaining segments $(O,A)$ and $(O,C)$. 
	On $(O,A)$, it was known that \eqref{eq:OS} holds when $\alpha < p$.
	Moreover, at the critical summability exponent $\alpha = \alpha_* = p$, the failure of \eqref{eq:OS restricted} for all $\beta>1$ was proved in \cite{Bez et al 2019}.
	More precisely, the failure when $d\ge 3$ was proved, but their proof also works when $d=2$.
	Our third main result claims that \eqref{eq:OS restricted} with $\beta=1$ holds on $(O,A)$.
	\begin{theorem}[On \text{$(O,A)$} when \text{$d \ge 2$}]\label{thm:OA}
		Let $d\geq2$ and $0<s<d/2$.  Set
		\begin{equation}\label{eq:p and q}
			p=\frac{d+1}{d-2s},
			\qquad
			q=\frac{d(d+1)}{(d-1)(d-2s)}, 
			\qquad 
			\alpha_* = p,
			\qquad
			\beta = 1.
		\end{equation}
		Then, the estimate \eqref{eq:OS restricted} holds for every $\nu\in\ell^{p,1}$ and for every orthonormal system $(f_j)_j$ in $\dot{H}^s_x(\R^d)$.
	\end{theorem}
	
	On $(O,C)$, Bez, Kinoshita, and Shiraki \cite{Bez et al 2024} obtained \eqref{eq:OS} if $\alpha < p$.
	When $d\ge2$, they also obtained the restricted weak-type estimate
	\cite[Section 3.1]{Bez et al 2024}
	\begin{equation}\label{eq:OC weak}
		\bigg\|\sum_j \nu_j |U(t)f_j|^2\bigg\|_
		{L_t^{p,\infty}(\R;L_x^\infty(\R^d))}
		\lesssim \|\nu\|_{\ell^{p,1}}.
	\end{equation}
	At the critical summability exponent $\alpha = \alpha_* = p$, it was proved that \eqref{eq:OS} fails in \cite{Gong 2026- 2}.
	The same counterexample proves that \eqref{eq:OS restricted} fails for all $\beta>1$ when $\alpha = p$. 
	Our fourth main result says that \eqref{eq:OS restricted} with $\beta=1$ holds on $(O,C)$.
	
	\begin{theorem}[On \text{$(O,C)$} when \text{$d \ge 2$}]\label{thm:OC}
		Let $d\ge2$, $p\in (1,\infty)$ and $q=\infty$.
		Set $s=d/2-1/p$, $\alpha_* = p$, and $\beta = 1$.
		Then, the estimate \eqref{eq:OS restricted} holds for every $\nu\in \ell^{p,1}$ and for every orthonormal system $(f_j)_j$ in $\dot{H}^s_x(\R^d)$.
	\end{theorem}

\subsection{Review of known results}\label{subsec:history}
In this section, we review known results. In summary, we confirm that Question \ref{question} has already been solved except 
\begin{itemize}
	\item $A$ when $d=1$;
	\item $(A,D)$, $(O,A)$, and $(O,C)$ when $d=2$;
	\item $(O,A)$ and $(O,C)$ when $d \ge 3$.
\end{itemize}

\subsubsection{On $(O,C)$ when $d=1$}
When $d=1$, the segments $(O,A)$ and $(O,C)$ coincide. On this segment, Bez, Kinoshita, and Shiraki \cite[Theorem 2.1]{Bez et al 2024} proved that \eqref{eq:OS} holds for $1\leq\alpha<p$. On the other hand, the discussion following that theorem shows that the critical estimate fails even with $\ell^{p,1}$ on the right-hand side and a weak time norm on the left-hand side.

\subsubsection{On $(O,B]$ when $d\ge 1$}
The case $p=\infty$, corresponding to the segment $(O,B]$, is already understood.
Lieb’s orthonormal Sobolev inequality \cite{Lieb 1983}, together with the unitarity of $U(t)$ on $\dot H^s_x(\R^d)$, yields \eqref{eq:OS restricted} for $0<s<d/2$, $q=\alpha_* = d/(d-2s)$, and $\beta=1$.
By \cite[Proposition 6.1]{Bez et al 2019}, this estimate fails if $\ell^{q,1}$ is replaced by $\ell^{q,\beta}$ for any $\beta>1$. Thus, on $(O,B)$, Question \ref{question} falls under case $(ii)$ with $\alpha_*=q$, and the optimal Lorentz exponent is $\beta=1$. At $B$, where $s=0$ and $q=1$, conservation of the $L^2$-norm gives \eqref{eq:naive} with the optimal exponent $\alpha=1$.

\subsubsection{On $(B,D)$ when $d\ge 1$}
In the pioneering work \cite{Frank et al 2014}, Frank, Lewin, Lieb, and Seiringer proved \eqref{eq:OS} with $s=0$ for exponents $p, q, \alpha\in[1, \infty]$ such that
	\begin{equation}
		\frac2p+\frac dq = d, \qquad 1\leq q\leq 1+\frac{2}{d},\qquad \alpha\leq\frac{2q}{q+1},
	\end{equation}
	which corresponds to the line segment $[B,F]$ in Figure \ref{fig:1}.
	Frank and Sabin \cite{Frank Sabin 2017} extended the previous result and proved that \eqref{eq:OS} with $s=0$ holds in the full range $[B,A)$. More precisely, they proved that \eqref{eq:OS} holds if
	\begin{equation}
		\frac2p+\frac dq = d, \qquad 1\leq q<\frac{d+1}{d-1},\qquad \alpha\leq \frac{2q}{q+1}.
	\end{equation}
	Moreover, in \cite{Frank et al 2014}, it was proved that \eqref{eq:OS} fails if $\alpha > 2q/(q+1)$ for all $1\le q < (d+1)/(d-1)$.
	It was also proved in \cite{Frank et al 2014} that \eqref{eq:OS} fails at the endpoint
	\[
	p=\frac{d+1}{d}, \qquad q=\frac{d+1}{d-1},\qquad \alpha = \frac{d+1}{d}.
	\]	
	The second and third authors \cite{Hadama Toyoshima 2026+} proved the restricted-type estimate at $A$ for $d\ge 2$:
	\begin{equation}\label{eq:}
		\bigg\|\sum_j \nu_j |U(t)f_j|^2\bigg\|_{L^{\frac{d+1}{d}}_t(\R;L^{\frac{d+1}{d-1}}_x(\R^d))} \lesssim \|\nu\|_{\ell^{\frac{d+1}{d},1}}.
	\end{equation}
	As a consequence of this result, it was proved that \eqref{eq:OS} holds on $(A,D]$ when $\alpha \le p$ and $d\ge 3$.
	The optimality of these estimates was proved in \cite{Frank et al 2014, Frank Sabin 2016}.
	
	\subsubsection{On $[D,C)$ when $d\ge 3$}
	When $d\ge3$, the segment $[D,C)$ corresponds to $p=1$ and $d/(d-2)\le q<\infty$.
	On this segment, \eqref{eq:naive} gives \eqref{eq:OS} with $\alpha=1$.
	The time-translation argument in \cite[Section 4.8]{Bez et al 2019} also applies with $p=1$ and shows that $\alpha\le1$ is necessary.
	
	\subsubsection{In $\interior OBA \cup \interior OADC$ when $d\ge 1$}
	Bez, Hong, Lee, Nakamura, and Sawano \cite{Bez et al 2019} proved that \eqref{eq:OS} holds in $\interior OBA$ for $1\le \alpha \le \alpha_*$, and \eqref{eq:OS} fails for $\alpha>\alpha_*$, where
	\[
	\frac1\alpha_* = \frac{1}{dp} + \frac{1}{q}.
	\]
	 Hence, Question \ref{question} is completely solved in $\interior OBA$.
	
	On the other hand, in \cite{Bez et al 2019}, it was proved that \eqref{eq:OS} holds in $\interior OADC$ for $1\le \alpha < p$, and \eqref{eq:OS} fails for $\alpha>p$.
	Recently, Feng, Lee, and Lee \cite{Feng et al 2026-} proved orthonormal Strichartz estimates with respect to fractal measures in the time variable.
	Specializing their results to Lebesgue time gives, for $d\ge3$, the critical summability estimates in $\interior OADC$ for $\alpha = p$.
	Note that their theorem contains both the admissible line estimates obtained in \cite{Hadama Toyoshima 2026+} and the $p=2$ result of \cite{Gong 2026-}.
	Moreover, in \cite{Gong 2026- 2}, the first author proved \eqref{eq:OS} with $\alpha=p$ in $\interior{OAC}$ in $d=2$.
	In conclusion, Question \ref{question} has already been solved in $\interior OBA\cup \interior OADC$.

    \subsection{Organization of this article}
    In Section \ref{sec:preliminary}, we introduce necessary notation and definitions. In particular, vector-valued Besov spaces and their interpolation inequality (Section \ref{subsec:Besov}) are very important in the proof of Theorems \ref{thm:OA} and \ref{thm:OC}.
    In Sections \ref{sec:d=1}, \ref{sec:d=2}, \ref{sec:OA}, and \ref{sec:OC}, we prove Theorems \ref{thm:d=1}, \ref{thm:d=2}, \ref{thm:OA}, and \ref{thm:OC}, respectively. Since the proof strategies are different for each theorem, we explain the idea of each proof at the beginning of each section, rather than giving a comprehensive summary in the introduction.

	\section{Preliminaries}\label{sec:preliminary}	
    In this section, we collect basic notation and tools that will be used throughout this paper.
	
	\subsection{Basic notation}\label{subsec:notation}
	Define $U(t)=e^{it\Delta_x}:=\mathcal{F}^{-1} e^{-it|\xi|^2} \mathcal{F}$, where the Fourier and inverse Fourier transforms are defined by
	\[
	\mathcal{F}[f](\xi):=\frac{1}{(2\pi)^{d/2}} \int_{\R^d} e^{-ix\cdot \xi} f(x)dx, 
	\qquad \mathcal{F}^{-1}[g](x):=\frac{1}{(2\pi)^{d/2}} \int_{\R^d} e^{ix\cdot \xi} g(\xi)d\xi.
	\]

   For $I\subset \Z$, $1\leq\alpha< \infty$, and $1\le \beta\leq\infty$, we define 
	\begin{empheq}[left={\|\nu\|_{\ell^{\alpha,\beta}(I)}:=\empheqlbrace}]{alignat=2}
		&\bigg(\sum_{j\ge 1} \frac1j (j^{1/\alpha}\nu_j^*)^\beta\bigg)^{1/\beta}&\quad&(\mathrm{if}\ 1\leq\beta<\infty) \\
		&\sup_{j\ge 1} j^{1/\alpha} \nu_j^*&\quad&(\mathrm{if}\ \beta=\infty),
	\end{empheq}
	where $(\nu_j^*)_{j\ge 1}$ is the rearrangement of $(|\nu_j|)_{j\in I}$ in decreasing order.
    Moreover, we define $\ell^{\infty,\infty} := \ell^\infty$.

    Let $H$ be a complex and separable Hilbert space.
    For every compact operator $A:H \to H$, define the Schatten--Lorentz quasi-norm by
	\[
	\|A\|_{\mathfrak{S}^{\alpha,\beta}(H)} := \|a\|_{\ell^{\alpha,\beta}},
	\]
	where $a=(a_j)_{j\ge 1}$ is the singular value sequence of $A$.
    Note that $\mathfrak S^\alpha(H) := \mathfrak{S}^{\alpha, \alpha}(H)$ is the standard Schatten-$\alpha$ class.
    Finally, we define $\mathfrak S^\infty(H) = \mathfrak S^{\infty,\infty}(H):= \mathcal B(H)$, where $\mathcal B(H)$ is the space of all bounded linear operators on $H$.

\subsection{Duality principle}
	We introduce the duality principle, which can be found in \cite{Bez et al 2019} and \cite{Frank Sabin 2017} for example.
	\begin{proposition}\label{prop:Frank Sabin duality}
		Let $p,q,r,s,\alpha,\beta\in [1,\infty]$ satisfy the following.
        \begin{itemize}
            \item If $p=\infty$, then $r=\infty$. Similarly, $s=\infty$ if $q=\infty$, and $\beta=\infty$ if $\alpha =\infty$;
            \item If $p=1$, then $r=1$. Similarly, $s=1$ if $q=1$, and $\beta=1$ if $\alpha=1$.
        \end{itemize}
        Then, 
		\begin{equation}
			\bigg\|\sum_j\nu_j|U(t)f_j|^2\bigg\|_{L_t^{p,r} L_x^{q,s}}\lesssim\|\nu\|_{\ell^{\alpha,\beta}}
		\end{equation}
		holds for every orthonormal system $(f_j)_j$ in $L^2_x$ and for every sequence $\nu=(\nu_j)_j\in\ell^{\alpha,\beta}$ if and only if
		\begin{equation}
			\|WUU^\ast W\|_{\mathfrak{S}^{\alpha',\beta'}(L^2_{t,x})}\lesssim\|W\|_{L_t^{2p',2r'}L_x^{2q',2s'}}^2.
		\end{equation}
	\end{proposition}

	\subsection{Littlewood--Paley decomposition}\label{subsec:LP}
	We can choose $\widehat\psi\in C_c^\infty(\R^d)$ such that $0\le \widehat\psi \le1$, 
    $\supp \widehat\psi \subset \{1/2\le |\xi|\le 2\}$, and 
	\[
	\sum_{k\in\Z}\widehat\psi_k(\xi)=1
	\qquad (\xi\in\R^d\setminus\{0\}),
	\]
    where $\widehat\psi_k(\xi):=\widehat\psi(2^{-k}\xi)$ for $k\in\Z$.
	For $f\in\mathcal S'(\R^d)$, define
	\[
	\widehat{P_k f}(\xi):=\widehat\psi_k (\xi)\widehat f(\xi).
	\]
	These operators satisfy
	\[
	\operatorname{supp}\widehat{P_k f}\subset\{\xi\in\R^d:2^{k-1}\le|\xi|\le2^{k+1}\},
	\qquad
	P_k P_{k'} f=0
	\]
    if $|k-k'|\ge 2$.
	The homogeneous Littlewood--Paley decomposition is
	\[
	f=\sum_{k\in\Z}P_k f,
	\]
	with convergence in $\mathcal S'(\R^d)$ modulo polynomials.
	If $f\in L^q_x$ with $1<q<\infty$, the series converges to $f$ in $L^q_x$.
	
	Define $\chi\in C_c^\infty(\R^d)$ by $\chi:= 1- \sum_{k\ge 1}\widehat\psi_k$.
	Let $P_{\le 0}:= \mathcal{F}^{-1} \chi \mathcal{F}$. 
	The inhomogeneous Littlewood--Paley decomposition is
	\[
	f=P_{\le 0}f+\sum_{k\ge 1}P_k f,
	\]
	with convergence in $\mathcal S'(\R^d)$ for every $f\in\mathcal S'(\R^d)$.
	If $f\in L^q_x$ with $1\le q<\infty$, this decomposition also converges in $L^q_x$.
	
	For $1<q<\infty$ and $f\in L^q_x$, both decompositions
	yield the square-function estimates
	\begin{align}
		\|f\|_{L^q_x}
		\sim \bigg\|\bigg(\sum_{k\in\Z}|P_k f|^2\bigg)^{1/2} \bigg\|_{L^q_x}
		\sim \bigg\| \bigg(|P_{\le 0} f|^2+\sum_{k\ge 1} |P_k f|^2\bigg)^{1/2} \bigg\|_{L^q_x}.
	\end{align}
    
    \subsection{Vector-valued Besov spaces and interpolation}\label{subsec:Besov}
    Note that we can write 
	\begin{equation}
		(P_k f)(x) = \frac{1}{(2\pi)^{d/2}}\int_{\R^d} 2^{dk} \psi(2^k(x-y))f(y)dy, 
	\end{equation}
	where $\psi \in \mathcal{S}(\R^d)$ is independent of $k\in \Z$.
	We write $\widetilde{P}_k = P_{k-1} + P_k + P_{k+1}$.
	
	Let $H$ be a Hilbert space.
	For every $f:\R^d \to H$ and for every Fourier multiplier $m(D)$, we define $m(D)f:\R^d\to H$ by
	\begin{equation}
		[m(D)f](x) := \frac{1}{(2\pi)^{d/2}}\int_{\R^d} \check{m}(x-y)f(y)dy \in H.
	\end{equation}
	We now introduce the $H$-valued homogeneous Besov space by
	\begin{equation}\label{eq:H valued Besov}
		\|f\|_{\dot B^s_{q,r;x}H}:=\big\|2^{ks} \|P_kf\|_{L_x^q H}\big\|_{\ell_k^r(\Z)}.
	\end{equation}
    The following interpolation inequality will be used in Sections \ref{sec:OA} and \ref{sec:OC}.

\begin{lemma}\label{lem:Besov interpolation}
Let $H$ be a Hilbert space.
Let $1<q_0\le q\le\infty$, $0\le s_0<s<s_1$, and $\theta \in (0,1)$ satisfy
\[
\frac1q=\frac{1-\theta}{q_0} + \frac{\theta}{\infty}, \qquad s=(1-\theta)s_0 + \theta s_1.
\]
Then, 
\begin{equation}\label{eq:general}
\||D|^{-s}f\|_{L^q_x H}
\lesssim \left\|\sup_{k\in\Z}2^{-s_0k}\|P_kf(x)\|_H\right\|_{L^{q_0}_x}^{1-\theta}
\|f\|_{\dot B_{\infty,\infty;x}^{-s_1}H}^\theta.
\end{equation}
In particular, when $s_0=0$,
\begin{equation}\label{eq:original}
\||D|^{-s}f\|_{L^q_xH}
\lesssim\|f\|_{L^{q_0}_xH}^{1-\theta}\|f\|_{\dot B_{\infty,\infty;x}^{-s_1}H}^\theta,
\qquad \theta=\frac{s}{s_1}.
\end{equation}
Moreover, when $q_0=\infty$,
\begin{equation}\label{eq:besov}
\||D|^{-s}f\|_{L^\infty_xH}
\lesssim\|f\|_{\dot B_{\infty,\infty;x}^{-s_0}H}^{1-\theta}\|f\|_{\dot B_{\infty,\infty;x}^{-s_1}H}^\theta.
\end{equation}
\end{lemma}

\begin{proof}
Let $M$ denote the scalar Hardy--Littlewood maximal operator.
Recall that $\widetilde P_k=P_{k-1}+P_k+P_{k+1}$ and $\widetilde P_kP_k=P_k$.
The integral kernel of $\widetilde P_k|D|^{-s}$ can be written as
\[
\widetilde K_k(x-y)=2^{k(d-s)}\widetilde K(2^k(x-y)),
\qquad \widetilde K\in\mathcal S(\R^d).
\]
First, the $L^1$ bound for this kernel gives
\begin{align*}
\|P_k|D|^{-s}f(x)\|_H
&=\|\widetilde P_k|D|^{-s}P_kf(x)\|_H
\le \|\widetilde K_k\|_{L^1}\|P_kf\|_{L_x^\infty H}\\
&\lesssim 2^{-ks}\|P_kf\|_{L_x^\infty H}
\le 2^{k(s_1-s)}\|f\|_{\dot B_{\infty,\infty;x}^{-s_1}H}.
\end{align*}
Consequently,
\begin{equation}\label{eq:low}
\sum_{k\le N}\|P_k|D|^{-s}f(x)\|_H
\lesssim2^{N(s_1-s)}\|f\|_{\dot B_{\infty,\infty;x}^{-s_1}H}.
\end{equation}
For the high-frequency terms, the same kernel gives
\begin{align*}
\|P_k|D|^{-s}f(x)\|_H
&\le 2^{-ks}\int_{\R^d}
2^{kd}|\widetilde K(2^k(x-y))|\|P_kf(y)\|_Hdy\\
&\le 2^{-k(s-s_0)}\int_{\R^d}
2^{kd}|\widetilde K(2^k(x-y))|
\sup_{j\in\Z}2^{-s_0j}\|P_jf(y)\|_Hdy\\
&\lesssim 2^{-k(s-s_0)}
M\bigg(\sup_{j\in\Z}2^{-s_0j}\|P_jf(\cdot)\|_H\bigg)(x),
\end{align*}
where we used \cite[Proposition, p.~24]{Stein 1993 textbook}.
Thus,
\begin{equation}\label{eq:high}
\sum_{k>N}\|P_k|D|^{-s}f(x)\|_H
\lesssim 2^{-N(s-s_0)}M\bigg(\sup_{j\in\Z}2^{-s_0j}\|P_jf(\cdot)\|_H\bigg)(x).
\end{equation}

Combining \eqref{eq:low} and \eqref{eq:high}, we obtain
\begin{equation}\label{eq:split}
\begin{aligned}
\bigl\||D|^{-s}f(x)\bigr\|_H
&\le\sum_{k\in\Z}\|P_k|D|^{-s}f(x)\|_H\\
&\lesssim 2^{N(s_1-s)}\|f\|_{\dot B_{\infty,\infty;x}^{-s_1}H}+2^{-N(s-s_0)}
M\bigg(\sup_{j\in\Z}2^{-s_0j}\|P_jf(\cdot)\|_H\bigg)(x).
\end{aligned}
\end{equation}
Optimizing with respect to $N$,
\begin{equation}\label{eq:pointwise}
\bigl\||D|^{-s}f(x)\bigr\|_H
\lesssim \bigg[M\bigg(\sup_{j\in\Z}2^{-s_0j}
\|P_jf(\cdot)\|_H\bigg)(x)\bigg]^{1-\theta}
\|f\|_{\dot B_{\infty,\infty;x}^{-s_1} H}^\theta.
\end{equation}
Taking the $L^q_x$ norm and using the $L^{q_0}_x$ boundedness of $M$, including $q_0=\infty$, gives
\begin{align*}
\||D|^{-s}f\|_{L^q_x H}
&\lesssim \|f\|_{\dot B_{\infty,\infty;x}^{-s_1} H}^\theta
\bigg\|M\bigg(\sup_{j\in\Z}2^{-s_0j}
\|P_jf(\cdot)\|_H\bigg)\bigg\|_{L^{q_0}_x}^{1-\theta}\\
&\lesssim \|f\|_{\dot B_{\infty,\infty;x}^{-s_1}H}^\theta
\bigg\|\sup_{j\in\Z}2^{-s_0j}\|P_jf(\cdot)\|_H\bigg\|_{L^{q_0}_x}^{1-\theta}.
\end{align*}
This proves \eqref{eq:general}.

When $s_0=0$, again by \cite[Proposition, p.~24]{Stein 1993 textbook}, 
\[
\sup_{k\in\Z}\|P_kf(x)\|_H
\lesssim M\bigl(\|f(\cdot)\|_H\bigr)(x).
\]
Hence
\[
\left\|\sup_{k\in\Z}\|P_kf(\cdot)\|_H\right\|_{L^{q_0}_x}
\lesssim \|f\|_{L^{q_0}_x H},
\]
which proves \eqref{eq:original}. Finally, when $q_0=\infty$,
\[
\left\|\sup_{k\in\Z}2^{-s_0k}\|P_kf(\cdot)\|_H\right\|_{L^\infty_x}
=\|f\|_{\dot B_{\infty,\infty;x}^{-s_0}H},
\]
and \eqref{eq:besov} follows.
\end{proof}

	\section{At \texorpdfstring{$A$}{A} when \texorpdfstring{$d=1$}{d=1}}\label{sec:d=1}
	In this section, we prove Theorem \ref{thm:d=1}.
	
	\subsection{Idea of the proof}
    We follow the strategy of Bez, Lee, and Nakamura \cite{Bez et al 2020}.
    In \cite{Bez et al 2020}, the authors obtain the weak-time estimate by applying the triangle inequality to $\sum_{l\ge1}W_1T_lW_2$. However, this does not use the fact that the kernel of $W_1T_lW_2$ is supported where $|t-\tau|\sim2^l$. They suggested using this separation to prove
    \[
    \bigg\|\sum_{l\ge1}W_1T_lW_2\bigg\|_{\mathfrak S^{b'}(L^2_{t,x})}
    \lesssim \bigg( \sum_{l\ge1}\|W_1T_lW_2\|_{\mathfrak S^{b'}(L^2_{t,x})}^{b}\bigg)^{1/b},
    \qquad 1\le b\le2.
    \]
They also checked $b=1,2$, but left open how to interpolate between them.

We prove that indeed we can interpolate them. Our main contribution is Lemma \ref{lem:interpolation}.
We keep the condition $|t-\tau|\sim2^l$ through smooth kernel cutoffs $M_l$ and consider
\[
(A_l)_l\longmapsto\sum_{l\ge1}M_l(A_l)
\]
for arbitrary $(A_l)_l$. The bounds at $b=1,2$ now apply to the same linear map and can be interpolated. Since $M_l(W_1T_lW_2)=W_1T_lW_2$, this gives the desired estimate. The key is to put the support condition into $M_l$, so that the inputs $A_l$ no longer have to share the special form $W_1T_lW_2$.
    	
	\subsection{Setup and review of the argument in \cite{Bez et al 2020}}
    Define the low-frequency cutoff $P$ such that
	\begin{equation*}
		\widehat{Pf}(\xi)=\chi(\xi)\widehat{f}(\xi),\qquad \text{where }\chi\in C_c^\infty(\mathbb{R}),\qquad 0\leq\chi\leq1,\qquad \chi(0)=1.
	\end{equation*}
	Since the estimate \eqref{eq:OS} is scaling invariant, it suffices to prove the localized estimate
	\begin{equation}\label{eq: frequency localized}
		\bigg\|\sum_j\nu_j\lvert U(t)Pf_j\rvert^2\bigg\|_{L_t^2 L_x^\infty}
		\lesssim_{\chi,\alpha}\|\nu\|_{\ell^\alpha}
	\end{equation}
	for every $\nu\in \ell^\alpha$, every orthonormal system $(f_j)_{j=1}^\infty$ in $L^2_x$, and for every $1\leq\alpha<2$.
	Indeed, we can obtain the desired estimate by \eqref{eq: frequency localized} with a standard rescaling argument. 
	
	Let $\alpha\in (4/3,2)$. Choose an annular cutoff $\eta\in C_c^\infty((-2,-1/2)\cup(1/2,2))$ and $\eta_0\in C_c^\infty(\R)$ such that\begin{equation*}
		\eta_0(t)+\sum_{l\geq1}\eta(2^{-l}t)=1.
	\end{equation*}
	Define $T$ by
	\[
	T=UP^2 U^*: f(t,x) \mapsto U(t)P\int_\R P U(\tau)^* f(\tau)d\tau.
	\]
	We decompose $T= T_0 + \sum_{l\ge 1} T_l$, where
	\[
	T_0(t,x;\tau,y) := \eta_0(t-\tau) T(t,x;\tau,y), \qquad T_l(t,x;\tau,y) := \eta(2^{-l}(t-\tau))T(t,x;\tau,y).
	\]
	Here, we identify an operator $A:L^2_{t,x} \to L^2_{t,x}$ with its integral kernel $A(t,x;\tau,y)$.
	By Proposition \ref{prop:Frank Sabin duality}, we see that the inequality \eqref{eq: frequency localized} is equivalent to
	\begin{equation*}
		\|W_1 T W_2\|_{\mathfrak{S}^{\alpha'}(L^2_{t,x})}
		=\bigg\|W_1\bigg(T_0+\sum_{l\geq1}T_l\bigg)W_2\bigg\|_{\mathfrak{S}^{\alpha'}(L^2_{t,x})}
		\lesssim\|W_1\|_{L_t^4L_x^2}\|W_2\|_{L_t^4L_x^2}.
	\end{equation*}
	In \cite[Section 3]{Bez et al 2020}, the authors provided a strong-type estimate
	\begin{equation}\label{eq:near-part}
		\|W_1T_0W_2\|_{\mathfrak{S}^2(L^2_{t,x})}\lesssim\|W_1\|_{L_t^{4}L_x^2}
		\|W_2\|_{L_t^{4}L_x^2}.
	\end{equation}
	A weak-time estimate
	\begin{equation}\label{eq:far part weak}
		\sum_{l\geq1}\|W_1T_lW_2\|_{\mathfrak{S}^{\alpha'}(L^2_{t,x})}
		\lesssim_\alpha
		\|W_1\|_{L_t^{4,2}L_x^2}
		\|W_2\|_{L_t^{4,2}L_x^2}
	\end{equation}
	is also deduced for $\alpha\in (4/3,2)$ (see \cite[Proof of Proposition 3.1]{Bez et al 2020}).
	Thus, the weak-time estimate
	\begin{equation*}
		\bigg\|\sum_j\nu_j\lvert U(t)Pf_j\rvert^2\bigg\|_{L_t^{2,\infty}L_x^\infty}
		\lesssim_{\chi,\alpha}\|\nu\|_{\ell^\alpha}
	\end{equation*}
	follows immediately by applying the triangle inequality.
	
	\subsection{Proof of Theorem \ref{thm:d=1}}
	We first establish the sequence estimate underlying \cite[(3.14)]{Bez et al 2020}.	
	\begin{lemma}\label{lem: bergh}
		For $4/3<b<2$, the bilinear transformation
	    \begin{equation}\label{eq: bergh lofstrom}
			\mathcal T(W_1,W_2)=(W_1T_lW_2)_{l\ge1},\qquad\mathcal T:
			L_t^{4,2b}L_x^2\times L_t^{4,2b}L_x^2
			\longrightarrow
			\ell^b(\mathfrak S^{b'}(L^2_{t,x}))
		\end{equation}is bounded.
	\end{lemma}
	\begin{proof}
		Fix $4/3<b<2$. We use the estimate \cite[(3.8)]{Bez et al 2020}
		\begin{equation*}
			\|W_1T_lW_2\|_{\mathfrak{S}^{b'}(L^2_{t,x})}
			\lesssim
			2^{\,l(1/2-1/q_1-1/q_2)}
			\|W_1\|_{L_t^{q_1}L_x^2}
			\|W_2\|_{L_t^{q_2}L_x^2},
		\end{equation*}
		for $(1/q_1,1/q_2)$ sufficiently close to $(1/4,1/4)$.
		Choose $\varepsilon>0$ sufficiently small and set
		\begin{equation*}
			a_0=\frac14+\varepsilon,\qquad
			a_1=\frac14-3\varepsilon,\qquad
			A_j=L_t^{1/a_j}L_x^2,\qquad j=0,1.
		\end{equation*}
		Define the weighted sequence norm by $\|(F_l)\|_{\ell^\infty_s(\mathfrak{S}^{b'}(L^2_{t,x}))}=\sup_{l\ge1}2^{sl}\|F_l\|_{\mathfrak{S}^{b'}(L^2_{t,x})}$. Put
		\begin{equation*}
			C_0=\ell^\infty_{2\varepsilon}(\mathfrak{S}^{b'}(L^2_{t,x})),
			\qquad
			C_1=\ell^\infty_{-2\varepsilon}(\mathfrak{S}^{b'}(L^2_{t,x})).
		\end{equation*}Then, we obtain the following boundedness:
		\begin{equation*}
			\mathcal T:A_0\times A_0\to C_0,
			\qquad
			\mathcal T:A_0\times A_1\to C_1,
			\qquad
			\mathcal T:A_1\times A_0\to C_1.
		\end{equation*}
		Using \cite[Exercise 5(b), Section 3.13]{Bergh Lofstrom 1976 textbook} yields a bounded bilinear transformation
		\begin{equation}\label{eq: bilinear bound}
			\mathcal T:
			(A_0,A_1)_{1/4,2b}
			\times
			(A_0,A_1)_{1/4,2b}
			\longrightarrow
			(C_0,C_1)_{1/2,b}.
		\end{equation}
		
		The standard real interpolation identity gives\begin{equation*}
			\bigl(L_t^{1/a_0}L_x^2,L_t^{1/a_1}L_x^2\bigr)_{\theta,s}
			=
			L_t^{q,s}L_x^2,
			\qquad
			\frac1q=(1-\theta)a_0+\theta a_1.
		\end{equation*}
		With $\theta=1/4$, we have $(A_0,A_1)_{1/4,2b}=L_t^{4,2b}L_x^2$. On the other hand, \cite[Theorem 5.6.1]{Bergh Lofstrom 1976 textbook} says
		\begin{equation*}
			\bigl(\ell^\infty_{s_0}(\mathfrak{S}^{b'}(L^2_{t,x})),\ell^\infty_{s_1}(\mathfrak{S}^{b'}(L^2_{t,x}))\bigr)_{\theta,r}
			=
			\ell^r_{(1-\theta)s_0+\theta s_1}(\mathfrak{S}^{b'}(L^2_{t,x})),
			\qquad s_0\ne s_1.
		\end{equation*}Therefore, we have $(C_0,C_1)_{1/2,b}
		=\ell^b_0(\mathfrak{S}^{b'})
		=\ell^b(\mathfrak{S}^{b'})$.
		In conclusion, \eqref{eq: bilinear bound} gives the boundedness in \eqref{eq: bergh lofstrom} as desired.
	\end{proof}
	
	We now focus on the following inequality:
	\begin{equation}\label{eq:pieces-lb}
		\bigg\|\sum_{l\geq1}W_1T_lW_2\bigg\|_{\mathfrak{S}^{b'}(L^2_{t,x})}
		\lesssim\bigg(\sum_{l\geq1}\|W_1T_lW_2\|_{\mathfrak{S}^{b'}(L^2_{t,x})}^b\bigg)^{1/b}
		\lesssim_b
		\|W_1\|_{L_t^{4,2b}L_x^2}
		\|W_2\|_{L_t^{4,2b}L_x^2}.
	\end{equation}
	The second inequality with $4/3<b<2$ follows directly from Lemma \ref{lem: bergh}.
	The first inequality with $4/3<b<2$ is the main difficulty identified in \cite{Bez et al 2020}.
	Let us first see why \eqref{eq:pieces-lb} gives \eqref{eq: frequency localized} as desired. 
	
	\begin{proof}[Proof of Theorem \ref{thm:d=1}, assuming \eqref{eq:pieces-lb}]
	Localizing the estimate in \cite[(3.12)]{Bez et al 2020} gives 
	\begin{equation}\label{eq: int 1}
		\bigg\|W_1\bigg(T_0+\sum_{l\geq1}T_l\bigg)W_2\bigg\|_{\mathfrak{S}^{4}(L^2_{t,x})}
		\lesssim\|W_1\|_{L_t^{4,8}L_x^2}\|W_2\|_{L_t^{4,8}L_x^2}.
	\end{equation}
	Now, suppose \eqref{eq:pieces-lb} holds for every $4/3<b<2$.
	Together with \eqref{eq:near-part}, we obtain
	\begin{equation}\label{eq: int 2}
		\bigg\|W_1\bigg(T_0+\sum_{l\geq1}T_l\bigg)W_2\bigg\|_{\mathfrak{S}^{b'}(L^2_{t,x})}
		\lesssim\|W_1\|_{L_t^{4,2b}L_x^2}\|W_2\|_{L_t^{4,2b}L_x^2}.
	\end{equation}
	Fix $\alpha\in(4/3,2)$, and choose
	\begin{equation*}
		b=\frac{8\alpha-8}{3\alpha-2},\qquad\theta=\frac{b}{4-b}.
	\end{equation*}
	We have $4/3<b<2$ and
	\begin{equation*}
		\frac{1}{2}=\frac{1-\theta}{4}+\frac{\theta}{b},\qquad\frac{1}{\alpha'}=\frac{1-\theta}{4}+\frac{\theta}{b'},\qquad 0<\theta=\frac{2}{\alpha'}<1.
	\end{equation*}
	Thus, complex interpolation gives $[\mathfrak{S}^4,\mathfrak{S}^{b'}]_\theta=\mathfrak{S}^{\alpha'}$. Moreover, the identity \begin{equation}\label{eq:special interpolation}
		\bigl[L_t^{4,8}(L_x^2),L_t^{4,2b}(L_x^2)\bigr]_\theta
		=L_t^4(L_x^2)
	\end{equation}
	is true for $1\leq b<2$, which, for example, can be deduced from \cite[(1.3)]{Hernandez Soria 1991}. In conclusion, we interpolate \eqref{eq: int 1}, \eqref{eq: int 2}, and thus obtain
	\begin{equation}\label{eq: strong near part}
		\bigg\|W_1\bigg(T_0+\sum_{l\geq1}T_l\bigg)W_2\bigg\|_{\mathfrak{S}^{\alpha'}(L^2_{t,x})}
		\lesssim
		\|W_1\|_{L_t^{4}L_x^2}
		\|W_2\|_{L_t^{4}L_x^2}
	\end{equation}
	for every $\alpha\in(4/3,2)$. This gives \eqref{eq: frequency localized}.	
	\end{proof}
	
	It remains to prove the first inequality of \eqref{eq:pieces-lb}.
	To prove it, we address the interpolation difficulty discussed after (3.13) in \cite{Bez et al 2020}.
	Choose another annular cutoff $\widetilde\eta\in C_c^\infty(\R)$ which satisfies
	\begin{equation*}
		\supp \widetilde\eta=\{t; 1/4\leq|t|\leq4\},\qquad \widetilde\eta\equiv1\text{ on }\supp\eta, \qquad 0\le \widetilde{\eta}\le 1.
	\end{equation*}
	Let $\widetilde\eta_l(t)=\widetilde\eta(2^{-l}t)$.
	For a bounded operator $A:L^2_{t,x} \to L^2_{t,x}$, define $\mathcal M_l(A)$ by
	\begin{equation}\label{eq: kernel}
		\mathcal M_l(A)(t,x;\tau,y) = \widetilde\eta_l(t-\tau)A(t,x;\tau,y).
	\end{equation}
	We also emphasize that
	\begin{equation*}
		\mathcal{M}_l(W_1T_lW_2)=W_1T_lW_2
	\end{equation*}
	for all $l\geq1$ since $\widetilde\eta = 1$ on $\supp \eta$.
	We finally prove the following lemma, which immediately implies the first inequality of \eqref{eq:pieces-lb}.

	\begin{lemma}\label{lem:interpolation}
		Let $1\leq b\leq2$. Then, for every sequence $\mathbf{A}=(A_l)_l \in\ell^b(\mathfrak{S}^{b'}(L^2_{t,x}))$, \begin{equation}\label{eq: int}
			\bigg\|\sum_{l\geq1}\mathcal{M}_l(A_l)\bigg\|_{\mathfrak{S}^{b'}(L^2_{t,x})}
			\lesssim_b\|\mathbf{A}\|_{\ell^b(\mathfrak{S}^{b'}(L^2_{t,x}))}.
		\end{equation}
	\end{lemma}
	
	\begin{proof}
		First, note that
		\begin{equation}\label{eq:bound of M}
			\|\mathcal{M}_l(A)\|_{\mathfrak S^r(L^2_{t,x})}
			\lesssim \|\widetilde \eta\|_{\mathcal{F} L^1} \|A\|_{\mathfrak{S}^r(L^2_{t,x})}.
		\end{equation}
		Indeed, since
		\[
		\mathcal{M}_l(A) = \frac1{\sqrt{2\pi}}\int_\R e^{it\omega} A e^{-it\omega} \mathcal{F}_{t\to \omega}[\widetilde{\eta}_l](\omega) d\omega,
		\]
		we obtain
		\begin{align}
		\|\mathcal{M}_l(A)\|_{\mathfrak{S}^r(L^2_{t,x})}
		&\lesssim  \int_\R \|A\|_{\mathfrak{S}^r(L^2_{t,x})} |\mathcal{F}_{t\to \omega}[\widetilde{\eta}_l](\omega)| d\omega \\
		&= \|\widetilde{\eta}_l\|_{\mathcal{F} L^1}\|A\|_{\mathfrak{S}^r(L^2_{t,x})}
		= \|\widetilde{\eta}\|_{\mathcal{F} L^1}\|A\|_{\mathfrak{S}^r(L^2_{t,x})}.
		\end{align}
		By the triangle inequality and \eqref{eq:bound of M}, 
		\begin{equation}\label{eq:int3}
		\begin{aligned}
			\bigg\|\sum_{l\ge 1}\mathcal{M}_l(A_l)\bigg\|_{\mathfrak{S}^\infty(L^2_{t,x})}
			&\leq\sum_{l\ge 1}\|\mathcal{M}_l(A_l)\|_{\mathfrak{S}^\infty(L^2_{t,x})}\\
			&\lesssim \|\widetilde\eta\|_{\mathcal{F} L^1}\sum_{l\ge 1}\|A_l\|_{\mathfrak{S}^\infty(L^2_{t,x})}
			\sim \|\mathbf{A}\|_{\ell^1(\mathfrak{S}^\infty(L^2_{t,x}))}.
		\end{aligned}
	\end{equation}
		
		Next, we consider the $\mathfrak S^2$-norm. We have
		\begin{equation}
			\bigg\|\sum_l\mathcal{M}_l(A_l)\bigg\|_{\mathfrak{S}^{2}(L^2_{t,x})}^2
			= \sum_{l,l'\ge 1} \iiiint_{\R^4} \mathcal{M}_l(A_l)(t,x;\tau,y)\overline{\mathcal{M}_{l'}(A_{l'})(t,x;\tau,y)} dtdxd\tau dy.
		\end{equation}
		Since each $\widetilde\eta_l$ is supported on the $l$-th temporal annulus, the terms on the right-hand side have uniformly finite overlap.
		Therefore, there exists a finite $L\in \mathbb N$ such that 
		\begin{equation}\label{eq:int4}
			\begin{aligned}
				\bigg\|\sum_l\mathcal{M}_l(A_l)\bigg\|_{\mathfrak{S}^{2}(L^2_{t,x})}^2
				&\le \sum_{|l-l'|\le L} \iiiint_{\R^4} |A_l(t,x;\tau,y) A_{l'}(t,x;\tau,y)| dtdxd\tau dy \\
				&\lesssim_L \|\mathbf A\|_{\ell^2(\mathfrak{S}^2(L^2_{t,x}))}^2.
			\end{aligned}
		\end{equation}

		Finally, for $1<b<2$, we take $0<\theta<1$ such that $1/b=1-\theta/2$. This gives $1/b'=\theta/2$ and
		\begin{equation*}			[\mathfrak{S}^\infty(L^2_{t,x}),\mathfrak{S}^2(L^2_{t,x})]_\theta=\mathfrak{S}^{b'}(L^2_{t,x}),
        \qquad[\ell^1(\mathfrak{S}^\infty(L^2_{t,x})),\ell^2(\mathfrak{S}^2(L^2_{t,x}))]_\theta=\ell^b(\mathfrak{S}^{b'}(L^2_{t,x})),
		\end{equation*}
		where we used \cite[Theorem 5.1.2]{Bergh Lofstrom 1976 textbook}.
		Using complex interpolation between \eqref{eq:int3} and \eqref{eq:int4}, we obtain \eqref{eq: int}.
	\end{proof}

	\section{On \texorpdfstring{$(A,D)$}{(A,D)} when \texorpdfstring{$d=2$}{d=2}}\label{sec:d=2}
	In this section, we show Theorem \ref{thm:d=2}.
	Throughout this section, $1<p<3/2$ and $q=p'$.
	
	\subsection{Idea of the proof}
    For $d\ge3$, we can obtain the estimates on $(A,D)$ by interpolating between $A$ and $D$ (see \cite{Hadama Toyoshima 2026+}). This does not work for $d=2$, since the Strichartz estimate for $(p,q)=(1,\infty)$ fails (see \cite{Montgomery-Smith 1998}). Thus, we directly apply the near-far decomposition argument, which was developed in \cite{Hadama Toyoshima 2026+}.
    
    In \cite{Hadama Toyoshima 2026+}, the time scale is chosen according to $V(t,x)\in L^q_t L^p_x$, and the energy estimate is used for short-time interactions. Our idea is to choose the time scale according to $\|V(t)\|_{L_x^p}$, using
    \[
    F_k=\bigl\{t:2^k\le\|V(t)\|_{L_x^p}<2^{k+1}\bigr\}.
    \]
    The local norm that enters the Strichartz estimate satisfies
    \[
    \bigl\|\one_{F_{k'}}\one_{[a,a+\ell(k,k'))}V\bigr\|_{L_t^qL_x^p}
    \le 2^{k'+1}\ell(k,k')^{1/q}.
    \]
    Thus, shortening $\ell(k,k')$ compensates for large values of $\|V(t)\|_{L_x^p}$. We choose $\ell(k,k')$ to obtain this local control while leaving enough decay in $|t-\tau|$ to sum the time-separated interactions. The condition $q>2p$ makes these requirements compatible.

    \subsection{A Schatten estimate with a time cutoff}

We identify integral operators with their operator-valued time kernels in the following statement.

\begin{lemma}\label{lem:dyadic}
For every $\chi,W_1,W_2$, let
\[
A_\chi(W_1,W_2) u(t):= W_1(t)\int_\R \chi(t-\tau)U(t-\tau)(W_2(\tau)u(\tau)) d\tau. 
\]
Then, 
\begin{equation}\label{eq:dyadic}
\|A_\chi(W_1,W_2)\|_{\mathfrak S^{2p}(L^2_{t,x})}
\lesssim \||\cdot|^{-1/p}\chi\|_{L_t^{\frac{2p}{2p-1}}}\|W_1\|_{L_t^{4p}L_x^{2p}}\|W_2\|_{L_t^{4p}L_x^{2p}}.
\end{equation}
\end{lemma}

   \begin{proof}
        In two dimensions, the Schr\"odinger kernel is
        \begin{equation}
            K_\tau(x,y)=\frac{1}{4\pi i\tau}\exp\left(\frac{i|x-y|^2}{4\tau}\right),
        \end{equation}
        so $|K_\tau(x,y)|\lesssim|\tau|^{-1}$.
        Therefore, by Young's inequality,
        \begin{align}
            \|A_\chi(W_1,W_2)\|_{\mathfrak S^2(L^2_{t,x})}^2
            &\lesssim\int_\R dt\int_{\R^2} dx \int_\R d\tau \int_{\R^2} dy \frac{|\chi(t-\tau)|^2|W_1(t,x)|^2|W_2(\tau,y)|^2}{|t-\tau|^2}\\
            &=\iint_{\R^2}\frac{|\chi (t-\tau)|^2\|W_1(t)\|_{L_x^2}^2\|W_2(\tau)\|_{L_x^2}^2}{|t-\tau|^2}dtd\tau\\
            &\lesssim \bigg\|\frac{\chi}{|\cdot|}\bigg\|_{L_t^2}^2\|W_1\|_{L_t^4L_x^2}^2\|W_2\|_{L_t^4L_x^2}^2.\label{eq:hs}
	   \end{align}
    
        On the other hand, the unitarity of the Schr\"odinger operator gives 
		\begin{equation}
			\|W_1(t)U(t-\tau)W_2(\tau)\|_{L_x^2\to L_x^2}\lesssim\|W_1(t)\|_{L_x^\infty}\|W_2(\tau)\|_{L_x^\infty},
		\end{equation}
		which yields the next estimate with Young's inequality:
		\begin{align}
			\|A_\chi(W_1,W_2)u\|_{L^2_{t,x}}
            &\lesssim \bigg\|\int|\chi (t-\tau)|\|W_1(t)\|_{L_x^\infty}\|W_2(\tau)\|_{L_x^\infty}\|u(\tau)\|_{L_x^2}d\tau\bigg\|_{L^2_t}\\
			&\leq\|W_1\|_{L_{t,x}^\infty}\|W_2\|_{L_{t,x}^\infty}\bigg\|\int|\chi (t-\tau)|\|u(\tau)\|_{L_x^2}d\tau\bigg\|_{L^2_t}\\
            &\le \|\chi\|_{L^1_t} \|W_1\|_{L_{t,x}^\infty}\|W_2\|_{L_{t,x}^\infty}\|u\|_{L_{t,x}^2}.
		\end{align}
        Therefore, 
        \begin{equation}\label{eq:op}
            \|A_\chi(W_1,W_2)\|_{\mathfrak{S}^\infty(L^2_{t,x})}
            \lesssim \|\chi\|_{L^1_t} \|W_1\|_{L_{t,x}^\infty}\|W_2\|_{L_{t,x}^\infty}.
        \end{equation}
        Then, we apply the trilinear complex interpolation with the weighted interpolation formula (see \cite[Theorem 5.5.3]{Bergh Lofstrom 1976 textbook})
        \[
        [L^1_t(\R;dt),L^2_t(\R;|t|^{-2}dt)]_{1/p}=L^{\frac{2p}{2p-1}}_t(\R;|t|^{-\frac2{2p-1}}dt)
        \]
        and obtain \eqref{eq:dyadic}.
    \end{proof}

    \subsection{Compactness}
    Choose an even nonnegative $\eta\in C_c^\infty(\R\setminus\{0\})$ with
    \begin{equation}
        \supp\eta\subset\{t:1/2\leq |t|\leq 2\},\qquad
        \eta_j(t)=\eta(2^{-j}t),\qquad
        \sum_{j\in\Z}\eta_j(t)=1\quad(t\ne0),
    \end{equation}
    and set
    \begin{equation}
        S_j g(t)=\int_\R\eta_j (t-\tau) U(t-\tau)g(\tau)d\tau.
    \end{equation}
    Let $V\ge0$ be compactly supported and bounded.
    Denote $W:=V^{1/2}$ and define
    \begin{equation}
        WUf(t,x):=W(t,x)U(t)f(x),\qquad
        (WU)^\ast f:=\int U(\tau)^\ast[W(\tau)f(\tau)]d\tau.
    \end{equation}
    We also denote $A:=WUU^*W$ and $U^*VU:=(WU)^*(WU)$.
       The following lemma is analogous to \cite[Lemma 3.6]{Hadama Toyoshima 2026+}.
    \begin{lemma}\label{lem:compactness}
        $A$ is positive and compact.
    \end{lemma}
    \begin{proof}
        Note that \eqref{eq:op} and \eqref{eq:hs} imply
        \begin{equation}
            \bigg\|\sum_{j\le J}WS_j W\bigg\|_{\mathfrak{S}^\infty(L^2_{t,x})}\lesssim2^J \|W\|_{L_{t,x}^\infty}^2,\qquad
            \bigg\|\sum_{j>J}WS_j W\bigg\|_{\mathfrak{S}^2(L^2_{t,x})}\lesssim2^{-J/2}\|W\|_{L_t^4L_x^2}^2.
        \end{equation}
        So $\sum_{j>J}WS_j W$ is compact, and thus
        \begin{equation}
            \bigg\|A-\sum_{j>J}WS_j W\bigg\|_{\mathfrak{S}^\infty(L^2_{t,x})}\lesssim2^J \|W\|_{L_{t,x}^\infty}^2\longrightarrow0\quad(J\to-\infty),
        \end{equation}
        which yields the compactness of $A$.
    \end{proof}

   \subsection{Near-far decomposition}
   Let $V\ge0$ be compactly supported, bounded, and satisfy
   \[
   \|V\|_{L_t^qL_x^p}=1
   \]
   with finitely many nonempty sets
    \begin{equation}
        F_k=\{t:2^k\le\|V(t)\|_{L_x^p}<2^{k+1}\},\qquad k\in\Z.
    \end{equation}
    
    \begin{remark}\label{rmk:W estimate}
        By definition,
        \begin{equation}\label{eq:levels-1}
            \|\one_{F_k}W\|_{L_t^{4p}L_x^{2p}}^2
            =\left(\int_{F_k}\|V(t)\|_{L_x^p}^{2p}dt\right)^{1/(2p)}\leq2^{k+1}|F_k|^{\frac1{2p}}.
        \end{equation}
        Also, since the sets $(F_k)_k$ are disjoint and $V$ is normalized,
        \begin{equation}\label{eq:levels-2}
            \sum_{k\in\Z}2^{qk}|F_k|\leq\sum_{k\in\Z}\int_{F_k}\|V(t)\|_{L_x^p}^qdt\leq\|V\|_{L_t^qL_x^p}^q=1.
        \end{equation}
    \end{remark}

    Recall that $A=W UU^* W$. We decompose $A$ into two parts: the near part and the far part.
    For $0<\varepsilon<1$ and $\lambda>0$, let
    \begin{equation}\label{eq:length}
        \ell(k,k')=(\varepsilon\lambda)^q 2^{-q(k+k')/2} 2^{q|k-k'|/4}.
    \end{equation}
    Note that
    \begin{equation}
    1=\bigg(1-\frac{|t-\tau|}{\ell(k,k')}\bigg)_+ +\min\bigg(1,\frac{|t-\tau|}{\ell(k,k')}\bigg),
    \end{equation}
    where $(\cdot)_+$ denotes the positive part.
    Define the near part $\mathscr{N}_{\varepsilon,\lambda}$ and the far part $\mathscr{F}_{\varepsilon,\lambda}$ by
    \begin{align}
        &\mathscr{N}_{\varepsilon,\lambda}(A):= \sum_{k,k'\in\Z}
        \bigg(1-\frac{|t-\tau|}{\ell(k,k')}\bigg)_+ \one_{F_k}(t)W(t)U(t-\tau)W(\tau)\one_{F_{k'}}(\tau),\label{eq:near-mask} \\
        &\mathscr{F}_{\varepsilon,\lambda}(A):= \sum_{k,k'\in\Z}
        \min\bigg(1,\frac{|t-\tau|}{\ell(k,k')}\bigg) \one_{F_k}(t)W(t)U(t-\tau)W(\tau)\one_{F_{k'}}(\tau).
    \end{align}

    We first estimate the far part as follows.
    \begin{lemma}[Far part]\label{lem:far}
        For $0<\varepsilon<1$ and $\lambda>0$,
        \begin{equation}\label{eq:far}
            \|\mathscr{F}_{\varepsilon,\lambda}(A)\|_{\mathfrak{S}^{2p}(L^2_{t,x})}\lesssim(\varepsilon\lambda)^{1-\frac{q}{2p}}.
        \end{equation}
    \end{lemma}
    \begin{proof}
    For each $k,k'\in \Z$, by Lemma \ref{lem:dyadic}, 
    \begin{align}
        \|\one_{F_k}\mathscr{F}_{\varepsilon,\lambda}(A)\one_{F_{k'}}\|_{\mathfrak{S}^{2p}(L^2_{t,x})}
        &\le \sum_{j\in\Z} \|A_{\chi_j}(\one_{F_k} W,\one_{F_{k'}} W) \|_{\mathfrak{S}^{2p}(L^2_{t,x})} \\
        &\lesssim \sum_{j\in\Z} \||\cdot|^{-1/p}\chi_j\|_{L_t^{\frac{2p}{2p-1}}} \|\one_{F_k}W\|_{L_t^{4p}L_x^{2p}}\|\one_{F_{k'}}W\|_{L_t^{4p}L_x^{2p}},
    \end{align}
    where
    \[
    \chi_j(t-\tau) = \eta_j(t-\tau)\min\bigg(1,\frac{|t-\tau|}{\ell(k,k')}\bigg).
    \]
    By the definition of $\chi_j$, \eqref{eq:length}, and Remark \ref{rmk:W estimate},
    \begin{align}
        \|\one_{F_k}\mathscr{F}_{\varepsilon,\lambda}(A)\one_{F_{k'}}\|_{\mathfrak{S}^{2p}(L^2_{t,x})}
        &\lesssim \sum_{j\in\Z} 2^{j(\frac1q-\frac1{2p})} \min\bigg(1,\frac{2^j}{\ell(k,k')}\bigg) 2^{(k+k')/2} |F_k|^{\frac{1}{4p}} |F_{k'}|^{\frac{1}{4p}} \\
        &\lesssim \ell(k,k')^{\frac1q-\frac1{2p}} 2^{(k+k')/2} |F_k|^{\frac{1}{4p}} |F_{k'}|^{\frac{1}{4p}} \\
        &=(\varepsilon\lambda)^{1-\frac{q}{2p}}2^{-\frac{(q-2p)}{8p}|k-k'|}(2^{qk}|F_k|)^{\frac{1}{4p}}(2^{qk'}|F_{k'}|)^{\frac{1}{4p}}.\label{eq:far block}
    \end{align}
    For each fixed $j$, the blocks indexed by $(k,k-j)$ have mutually orthogonal input spaces and mutually orthogonal output spaces. 
        Thus,
        \begin{align*}
            &\bigg\|\sum_k\one_{F_k}\mathscr{F}_{\varepsilon,\lambda}(A)\one_{F_{k-j}}\bigg\|_{\mathfrak{S}^{2p}(L^2_{t,x})}^{2p}
            =\sum_k\|\one_{F_k}\mathscr{F}_{\varepsilon,\lambda}(A)\one_{F_{k-j}}\|_{\mathfrak{S}^{2p}(L^2_{t,x})}^{2p}\\
            &\qquad \lesssim(\varepsilon\lambda)^{2p-q}2^{-\frac{q-2p}{4}|j|}\sum_k(2^{qk}|F_k|)^{1/2}(2^{q(k-j)}|F_{k-j}|)^{1/2}\\
            &\qquad \lesssim (\varepsilon\lambda)^{2p-q}2^{-\frac{q-2p}{4}|j|}\bigg(\sum_k2^{qk}|F_k|\bigg)^{1/2}\bigg(\sum_k2^{q(k-j)}|F_{k-j}|\bigg)^{1/2}\\
            &\qquad \lesssim (\varepsilon\lambda)^{2p-q}2^{-\frac{q-2p}{4}|j|}.
        \end{align*}
        For the first identity of the above estimates, see \cite[Introduction]{Carbery 2009}.
        Since $q>2p$, the triangle inequality gives
        \begin{align}
            \|\mathscr{F}_{\varepsilon,\lambda}(A)\|_{\mathfrak{S}^{2p}(L^2_{t,x})}
            \lesssim(\varepsilon\lambda)^{1-\frac{q}{2p}}\sum_{j\in\Z}2^{-\frac{(q-2p)}{8p}|j|}
            \lesssim(\varepsilon\lambda)^{1-\frac{q}{2p}}.
        \end{align}
    \end{proof}
    Next, we estimate the near part.
    \begin{lemma}[Low spectral part of the near part]
        For $0<\varepsilon<1$ and $\lambda>0$,
        let
        \[
        A_{\textup{low}}=WU\one_{[0,\lambda/\sqrt{\varepsilon}]}(U^\ast VU)U^\ast W.
        \]
        Then, we have
        \begin{equation}\label{eq:near low}
            \|\mathscr{N}_{\varepsilon,\lambda}(A_{\textup{low}})\|_{\mathfrak{S}^\infty(L^2_{t,x})}\lesssim\lambda\varepsilon^{1/4}.
        \end{equation}
    \end{lemma}
    \begin{proof}
        The spectral theorem gives
        \begin{equation}\label{eq:spectral-low}
            \|WU\one_{[0,\lambda/\sqrt\varepsilon]}(U^\ast VU)\|_{L_x^2\to L_{t,x}^2}^{2}
            \leq\frac{\lambda}{\sqrt\varepsilon}.
        \end{equation}
        The single-function Strichartz estimate for the non-endpoint admissible pair $(2p,2q)$, followed by H\"older's inequality and duality, gives
        \begin{multline}
            \|U^*W\one_{F_{k'}}\one_{[a,a+\ell(k,k'))}\|_{L_{t,x}^2\to L_x^2}^{2}
            \lesssim \|W\one_{F_{k'}}\one_{[a,a+\ell(k,k'))}\|_{L_t^{2q}L_x^{2p}}^2\\
            =\bigg(\int_{F_{k'}\cap[a,a+\ell(k,k'))}\|V(t)\|_{L_x^p}^{q}\,dt\bigg)^{1/q}
            \leq2^{k'+1}\ell(k,k')^{1/q}.\label{eq:local}
        \end{multline}
        The triangular cutoff is an average of interval cutoffs:
        \begin{equation}\label{eq:average}
            \bigg(1-\frac{|t-\tau|}{\ell(k,k')}\bigg)_+
            =\frac{1}{\ell(k,k')}\int_\R\one_{[a,a+\ell(k,k'))}(t)\one_{[a,a+\ell(k,k'))}(\tau)da.
        \end{equation}
        For each fixed $t$, the integrand is supported in an $a$-interval of length $\ell(k,k')$. 
        Cauchy--Schwarz in $a$, \eqref{eq:spectral-low}, and \eqref{eq:local} therefore imply, for $k\ge k'$,
        \begin{align*}
            \left\|\one_{F_k}\mathscr{N}_{\varepsilon,\lambda}(A_\textup{low})\one_{F_{k'}}g\right\|_{L_{t,x}^2}^{2}
            &\leq\frac1{\ell(k,k')}\int_\R\left\|\one_{F_k}\one_{[a,a+\ell(k,k'))}A_\textup{low}\one_{F_{k'}}\one_{[a,a+\ell(k,k'))}g\right\|_{L_{t,x}^2}^{2}da\\
            &\lesssim\frac{\lambda}{\sqrt\varepsilon}2^{k'}\ell(k,k')^{1/q-1}\int_\R\|\one_{[a,a+\ell(k,k'))}g\|_{L_{t,x}^2}^{2}da\\
            &=\frac{\lambda}{\sqrt\varepsilon}2^{k'}\ell(k,k')^{1/q}\|g\|_{L_{t,x}^2}^{2}.
        \end{align*}
        Moreover,
        \begin{equation}
            \frac{\lambda}{\sqrt\varepsilon}\,2^{k'}\ell(k,k')^{1/q}
            =\lambda^2\varepsilon^{1/2}2^{k'-(k+k')/2+|k-k'|/4}
            =\lambda^2\varepsilon^{1/2}2^{-|k-k'|/4}.
        \end{equation}
        For $k<k'$, use self-adjointness. 
        The same block orthogonality as above, now in operator norm, yields
        \begin{align}
            \|\mathscr{N}_{\varepsilon,\lambda}(A_\textup{low})\|_{\mathfrak{S}^\infty(L^2_{t,x})}
            &\leq\sum_{l\in\Z}\sup_{k\in\Z}\big\|\one_{F_k}\mathscr{N}_{\varepsilon,\lambda}(A_\textup{low})\one_{F_{k-l}}\big\|_{\mathfrak{S}^\infty(L^2_{t,x})}\\
            &\lesssim\lambda\varepsilon^{1/4}\sum_{l\in\Z}2^{-|l|/8}
            \lesssim\lambda\varepsilon^{1/4}.
        \end{align}
    \end{proof}

    \begin{lemma}[High spectral part of the near part]
    Let $\lambda_0>0$.
    The compactness of $A$ allows us to choose $\lambda_*\ge\lambda_0$ such that
    \begin{equation}\label{eq:max}
        \lambda_*^q\mathrm{rank}\one_{[\lambda_*,\infty)}(A)
        =\max_{\zeta\ge\lambda_0}\zeta^q\mathrm{rank}\one_{[\zeta,\infty)}(A).
    \end{equation}
    For $0<\varepsilon<1$, let
    \[
        A_{\textup{high}}=WU\one_{(\lambda_*/\sqrt{\varepsilon},\infty)}(U^\ast VU)U^\ast W.
    \]
    Then, we have
    \begin{equation}\label{eq:near high}
        \|\mathscr{N}_{\varepsilon,\lambda_*}(A_{\textup{high}})\|_{\mathfrak{S}^2(L^2_{t,x})}
        \lesssim \varepsilon^{q/4-1/2}\lambda_*\big(\rank\one_{[\lambda_*,\infty)}(A)\big)^{1/2}.
    \end{equation}
    \end{lemma}

    \begin{proof}
        By the polar decomposition,
        \begin{equation}
            A_\textup{high}=A\one_{(\lambda_*/\sqrt\varepsilon,\infty)}(A).
        \end{equation}
        Since the Hilbert--Schmidt norm is the $L^2$ norm of the integral kernel, we have
        \[
        \|\mathscr{N}_{\varepsilon,\lambda_*}(A_{\textup{high}})\|_{\mathfrak S^2(L^2_{t,x})}
        \le \|A_{\textup{high}}\|_{\mathfrak S^2(L^2_{t,x})}.
        \]
        List the (positive) eigenvalues $\zeta_j$ of $A_{\textup{high}}$ in nonincreasing order, with multiplicity.
        Since $\zeta_j>\lambda_*/\sqrt{\varepsilon}\ge\lambda_0$, \eqref{eq:max} gives
        \[
        j\zeta_j^q
        \le \zeta_j^q\rank\one_{[\zeta_j,\infty)}(A)
        \le \lambda_*^q\rank\one_{[\lambda_*,\infty)}(A).
        \]
        Since $q>2$, we obtain
        \[
        \begin{aligned}
        \|A_{\textup{high}}\|_{\mathfrak S^2(L^2_{t,x})}^2
        &=\sum_{j=1}^{\rank A_{\textup{high}}}\zeta_j^2
        \le \lambda_*^2 \bigl(\rank\one_{[\lambda_*,\infty)}(A)\bigr)^{2/q} \sum_{j=1}^{\rank A_{\textup{high}}} \frac{1}{j^{2/q}}\\
        &\lesssim \lambda_*^2\bigl(\rank\one_{[\lambda_*,\infty)}(A)\bigr)^{2/q}\bigl(\rank A_{\textup{high}}\bigr)^{1-2/q}.
        \end{aligned}
        \]
        Applying \eqref{eq:max} at $\lambda_*/\sqrt{\varepsilon}$ also gives
        \[
        \rank A_{\textup{high}}
        \le \rank\one_{[\lambda_*/\sqrt{\varepsilon},\infty)}(A)
        \le \varepsilon^{q/2}\rank\one_{[\lambda_*,\infty)}(A).
        \]
        Hence,
        \begin{equation}
        \|\mathscr{N}_{\varepsilon,\lambda_*}(A_{\textup{high}})\|_{\mathfrak S^2(L^2_{t,x})}^2
        \lesssim \varepsilon^{q/2-1}\lambda_*^2\rank\one_{[\lambda_*,\infty)}(A).
        \end{equation}
    \end{proof}

    Now, we prove a weak-type Schatten bound.
    \begin{proposition}\label{prop:weak}
       We have
        \begin{equation}\label{eq:weak}
            \|U^\ast VU\|_{\mathfrak{S}^{q,\infty}(L^2_x)}\lesssim\|V\|_{L_t^q L_x^p}.
        \end{equation}
    \end{proposition}
    \begin{proof}
        Using
        \begin{align*}
            \one_{[\lambda_*,\infty)}(A)A=\one_{[\lambda_*,\infty)}(A)\mathscr{F}_{\varepsilon,\lambda_*}(A)
            +\one_{[\lambda_*,\infty)}(A)\mathscr{N}_{\varepsilon,\lambda_*}(A_\textup{low})
            +\one_{[\lambda_*,\infty)}(A)\mathscr{N}_{\varepsilon,\lambda_*}(A_\textup{high}),
        \end{align*}
        we take the trace to obtain
        \begin{align*}
            \lambda_*&\mathrm{rank}\one_{[\lambda_*,\infty)}(A)\\
            &\leq\|\mathscr{F}_{\varepsilon,\lambda_*}(A)\|_{\mathfrak{S}^{2p}(L^2_{t,x})}\big(\mathrm{rank}\one_{[\lambda_*,\infty)}(A)\big)^{1-\frac{1}{2p}}+C_{p,q}\big(\varepsilon^{1/4}+\varepsilon^{q/4-1/2}\big)\lambda_*\mathrm{rank}\one_{[\lambda_*,\infty)}(A).
        \end{align*}
        Choose $\varepsilon=\varepsilon(p,q)>0$ so that the coefficient in parentheses is at most $1/2$. 
        After absorption, \eqref{eq:far} gives
        \begin{align*}
            \lambda_*^q\mathrm{rank}\one_{[\lambda_*,\infty)}(A)
            \lesssim \varepsilon^{2p-q}.
        \end{align*}
        By \eqref{eq:max} and the arbitrariness of $\lambda_0$, $\|A\|_{\mathfrak{S}^{q,\infty}(L^2_{t,x})}\lesssim1$.
        Since the nonzero eigenvalues of $A$ and $U^*VU$ coincide, we complete the proof.
    \end{proof}

    \subsection{Duality and interpolation}
    Now we prove Theorem \ref{thm:d=2}.
    Let $(f_j)_j$ be orthonormal in $L^2_x$.
    By the duality principle (Proposition \ref{prop:Frank Sabin duality}) and Proposition \ref{prop:weak},
    \begin{equation}\label{eq:restricted}
        \bigg\|\sum_j\nu_j|U(t)f_j|^2\bigg\|_{L_t^pL_x^q}\lesssim\|\nu\|_{\ell^{p,1}}.
    \end{equation}
    Choose $1<p_0<p<p_1<3/2$, $q_i=p_i'$, and $0<\theta<1$ with
    \begin{equation}
        \frac{1}{p}=\frac{1-\theta}{p_0}+\frac\theta{p_1},\qquad
        \frac{1}{q}=\frac{1-\theta}{q_0}+\frac\theta{q_1}.
    \end{equation}
    Then, we have \eqref{eq:restricted} for $(p,q)=(p_0,q_0), (p_1,q_1)$.
    We have
    \begin{align*}
        (\ell^{p_0,1},\ell^{p_1,1})_{\theta,p}&=\ell^p,\\
        (L_t^{p_0}L_x^{q_0},L_t^{p_1}L_x^{q_1})_{\theta,p}
        &=L_t^p\bigl((L_x^{q_0},L_x^{q_1})_{\theta,p}\bigr)
        =L_t^pL_x^{q,p}\hookrightarrow L_t^pL_x^q,
    \end{align*}
    where we used the Lions--Peetre theorem \cite{Lions Peetre 1964} (see also \cite{Cwikel 1974} and \cite[Proposition 5.1]{Lindemulder Lorist 2022}) for the second interpolation.
    By real interpolation, we conclude that
    \begin{equation}
        \bigg\|\sum_j\nu_j|U(t)f_j|^2\bigg\|_{L_t^pL_x^q}\lesssim\|\nu\|_{\ell^p}.
    \end{equation}

	\section{On \texorpdfstring{$(O,A)$}{(O,A)} when \texorpdfstring{$d\ge 2$}{d≧2}}\label{sec:OA}
	In this section, we prove Theorem \ref{thm:OA}.
	
    \subsection{Idea of the proof}\label{subsec:OA}
    The proof is given by interpolation between the endpoint estimate at $A$ given by \cite{Hadama Toyoshima 2026+} and the ``Sobolev embedding'' at $O$, namely, Lemma \ref{lem:O}.
    Note that, by Lieb's orthonormal Sobolev inequality \cite{Lieb 1983}, we have a useful bound
    \begin{equation}\label{eq:OS Sobolev}
    \bigg\|\sum_j \nu_j |f_j|^2 \bigg\|_{L^q_x}
    \lesssim \|\nu\|_{\ell^{q,1}}
    \end{equation}
    for all $1\le q < \infty$ and for every orthonormal system $(f_j)_j$ in $\dot H^{\frac{d}{2q'}}_x$.
    However, \eqref{eq:OS Sobolev} fails for $q=\infty$. 
    Moreover, it seems that we should not expect naively weakened versions of \eqref{eq:OS Sobolev}, for example,
    \begin{equation}
    \bigg\|\sum_j \nu_j |f_j|^2 \bigg\|_X
    \lesssim \|\nu\|_{\ell^\infty}
    \end{equation}
    with $X = \dot B^0_{\infty,\infty}$ or $X = \mathrm{BMO}$.

    Note that, for $\nu_j \ge 0$,
    \begin{equation}
        \bigg\| \sum_j \nu_j |U(t)f_j|^2 \bigg\|_{L^p_t L^q_x} \lesssim \|\nu\|_{\ell^{p,1}}
        \Longleftrightarrow \|\mu_j U(t) f_j\|_{L^{2p}_t L^{2q}_x \ell_j^2} \lesssim \|\mu_j\|_{\ell_j^{2p,2}},
    \end{equation}
    where $\mu_j = \sqrt{\nu_j}$.
    The key observation in this section is that we can prove an alternative bound for \eqref{eq:OS Sobolev} with $q=\infty$ if we use this form. Namely, for every orthonormal system $(f_j)_j$ in $L^2_x$ and for all $\mu_j \ge 0$, 
    \[
    \|\mu_j U(t) f_j\|_{L^\infty_t \dot B^{-d/2}_{\infty,\infty;x} \ell_j^2}\lesssim \sup_j \mu_j.
    \]
    See Section \ref{subsec:Besov} for vector-valued Besov spaces.
    Finally, we interpolate the above estimate and the endpoint estimate at $A$ by Lemma \ref{lem:Besov interpolation} to obtain the desired estimate.

    At the technical level, we first reduce to the problem of finite families with all coefficients equal to one (Lemma \ref{lem:coefficients}), which will also be used in Section \ref{sec:OC}.

	\subsection{Reduction}
	The following lemma reduces the proof of Theorem \ref{thm:OA} to a simpler estimate.
	\begin{lemma}\label{lem:coefficients}
		Let $p\in(1,\infty)$ and $q\in [1,\infty]$ and fix a family of functions $(u_j(t))_{j\ge 1}$.
		Suppose that
		\begin{equation}\label{eq:reduction}
			\|u_j(t)\|_{L_t^{2p}L_x^{2q}\ell^2_j(I)}
			\lesssim |I|^{\frac{1}{2p}}    
		\end{equation}
		for every finite $I\subset\N$.
		Then, for every $\nu\in \ell^{p,1}$,
		\[
		\bigg\|\sum_{j\ge 1}\nu_j|u_j(t)|^2\bigg\|_{L_t^pL_x^q}
		\lesssim_p \|\nu\|_{\ell^{p,1}}.
		\]
	\end{lemma}
	
	\begin{proof}
		We can assume that $\nu_j\ge \nu_{j+1}\ge 0$.
		Let $I_m := \{2^{m}, \dots, 2^{m+1}-1\}$.
		Then, we have 
		\begin{align*}
			\bigg\|\sum_{j\ge 1}\nu_j|u_j(t)|^2\bigg\|_{L_t^pL_x^q}
			&\leq \sum_{m\ge 0}\bigg\|\sum_{j\in I_m}\nu_j|u_j(t)|^2\bigg\|_{L_t^pL_x^q}\\
			&\lesssim \sum_{m\ge 0} \nu_{2^m} |I_m|^{1/p} = \sum_{m\ge 0} \nu_{2^m} 2^{m/p}.
		\end{align*}
		For $m\geq1$, monotonicity of $(\nu_j)$ gives
		\[
		\sum_{j=2^{m-1}}^{2^m-1}\nu_jj^{1/p-1} \gtrsim_p\nu_{2^m}2^{m/p}.
		\]
		Hence,
		\begin{align*}
			\sum_{m\ge 0} \nu_{2^m} 2^{m/p} \lesssim_p \sum_{j\ge 1} \nu_j j^{1/p-1} = \|\nu\|_{\ell^{p,1}}. 
		\end{align*}
	\end{proof}

	\subsection{Estimates at two endpoints}
	Let
	\begin{equation}\label{eq:A exponents}
		p_A:=\frac{d+1}{d},
		\qquad
		q_A:=\frac{d+1}{d-1}.
	\end{equation}
	
	The following is \cite[Corollary 1.8 $(ii)$]{Hadama Toyoshima 2026+}.
	\begin{lemma}[Estimate at \texorpdfstring{$A$}{A}]\label{lem:HT}
		Let $(f_j)_j$ be an orthonormal system in $L^2_x$,
		and let $(\nu_j)$ be a finitely supported scalar sequence.  Then
		\begin{equation}\label{eq:HT}
			\bigg\|\sum_j\nu_j|U(t)f_j|^2\bigg\|_{L_t^{p_A} L_x^{q_A}}
			\lesssim \|\nu\|_{\ell^{p_A,1}}.
		\end{equation}
	\end{lemma}
	
	We next prove the following lemma. We can regard it as the estimate at $O$.
	\begin{lemma}[Estimate at \text{$O$}]\label{lem:O}
		Let $(f_j)_{j\ge 1}$ be an orthonormal family in $L^2_x$, and let $\mu_j\geq0$.
		Then,
		\begin{equation}\label{eq:O}
			\|\mu_jf_j\|_{\dot{B}^{-d/2}_{\infty,\infty;x}\ell_j^2}  \lesssim \sup_{j\ge 1}\mu_j.
		\end{equation}
	\end{lemma}
	
	\begin{proof}
		Fix $k\in\Z$, and $x\in\R^d$.
		Bernstein's inequality and orthonormality give
		\begin{align*}
			\|[P_k \mu_jf_j](x)\|_{\ell_j^2}
			&=\sup_{\|c\|_{\ell^2}=1} \bigg|P_k\sum_{j\ge 1} c_j \mu_j f_j(x)\bigg|
			\lesssim 2^{kd/2} \sup_{\|c\|_{\ell^2}=1} \bigg\|\sum_{j\ge 1} c_j \mu_jf_j\bigg\|_{L^2_x} \\
			&= 2^{kd/2} \sup_{\|c\|_{\ell^2}=1} \bigg(\sum_{j\ge 1}\mu_j^2|c_j|^2\bigg)^{1/2}
			= 2^{kd/2}\sup_{j\ge 1}\mu_j.
		\end{align*}
		Multiplying by $2^{-kd/2}$ and taking the supremum over $x$ and $k$ proves \eqref{eq:O}.
	\end{proof}
	
	\subsection{Proof of Theorem \ref{thm:OA}}
	By Lemma \ref{lem:coefficients}, it suffices to prove \eqref{eq:reduction}.   
	Since $(f_j)_j$ is orthonormal in $\dot{H}^s_x$, $(|D|^sf_j)_j$ is orthonormal in $L^2_x$.
	Fix $I\subset \N$.
	Lemma \ref{lem:HT} gives
	\begin{align}
		\|U(t)|D|^s f_j\|_{L_t^{2p_A}L_x^{2q_A}\ell_j^2(I)}^2
		=\bigg\|\sum_{j\in I} |U(t)|D|^sf_j|^2\bigg\|_{L_t^{p_A}L_x^{q_A}}
		\lesssim |I|^{1/p_A}.
	\end{align}
	Lemma \ref{lem:O} gives the other endpoint directly:
	\begin{equation}\label{eq:O Fm}
		\|U(t)|D|^s f_j\|_{L_t^\infty \dot B^{-d/2}_{\infty,\infty;x} \ell^2_j(I)}
		\lesssim 1.
	\end{equation}
	Set
	\begin{equation}\label{eq:theta application}
		\theta:=\frac{s}{d/2}.
	\end{equation}
	Equations \eqref{eq:p and q} and \eqref{eq:A exponents} give
	\begin{equation}\label{eq:interpolated exponents}
		2p=\frac{2p_A}{1-\theta}, \qquad 2q=\frac{2q_A}{1-\theta}.
	\end{equation}
	Applying Lemma \ref{lem:Besov interpolation} to the spatial function $(U(t)|D|^sf_j)(x)$ for every fixed $t\in \R$,
	\begin{align}\label{eq:spatial-application}
		\|U(t)f_j\|_{L_x^{2q}\ell_j^2(I)}
		&= \||D|^{-s} U(t)|D|^sf_j\|_{L_x^{2q}\ell_j^2(I)}  \\
		&\lesssim \|U(t)|D|^sf_j\|_{L_x^{2q_A}\ell_j^2(I)}^{1-\theta}
		\|U(t)|D|^sf_j\|_{\dot{B}^{-d/2}_{\infty,\infty;x}\ell_j^2(I)}^{\theta} \\
		&\lesssim \|U(t)|D|^sf_j\|_{L_x^{2q_A}\ell_j^2(I)}^{1-\theta}.
	\end{align}
	Taking the $L_t^{2p}$ norm, 
	\begin{align}
		\|U(t)f_j\|_{L_t^{2p}L_x^{2q}\ell_j^2(I)}
		&\lesssim \|U(t)|D|^sf_j\|_{L^{2p_A}_t L_x^{2q_A}\ell_j^2(I)}^{1-\theta} 
		\lesssim |I|^{\frac{1}{2p}},
	\end{align}
	which completes the proof.

	\section{On \texorpdfstring{$(O,C)$}{(O,C)} when \texorpdfstring{$d\ge 2$}{d≧2}}\label{sec:OC}
	In this section, we prove Theorem \ref{thm:OC}.

    \subsection{Idea of the proof}
    Basically, we interpolate the ``Sobolev embedding'' at $O$ (namely, Lemma \ref{lem:O}) and the boundary Strichartz estimate at $C$ (see also Section \ref{subsec:OA}).
    At the technical level, we first reduce to the problem of finite families with all coefficients equal to one by Lemma \ref{lem:coefficients}.  
    Recall that the estimate at $O$ is given as an estimate for $\mu_j U(t)f_j$.
    Therefore, we should give an estimate at $C$ in the same form to use an interpolation argument.

    When $d\ge 3$, we indeed obtain a proper endpoint estimate at $C$, namely, Lemma \ref{lem:C}.
    This is a simple application of the endpoint Strichartz estimates by Keel and Tao \cite{Keel Tao 1998}.
    Hence, similarly to the proof of Theorem \ref{thm:OA}, the interpolation inequality (Lemma \ref{lem:Besov interpolation}) implies  Theorem \ref{thm:OC}.

    When $d =2$, we cannot expect Lemma \ref{lem:C} because the endpoint Strichartz estimate fails when $d=2$ (see \cite{Montgomery-Smith 1998, Tao 2006}). Hence, we interpolate estimates at $O$ and $(0,1-1/n)$ with $n\ge 2$.
    The argument for $d=2$ is inspired by \cite[Section 4.1]{Frank et al 2014}, \cite[Section 3]{Bez et al 2025}, and \cite[Section 6]{Gong 2026- 2}.

	\subsection{When \texorpdfstring{$d\ge 3$}{d≧3}}
	\begin{lemma}[Estimate at $C$]\label{lem:C}
		Let $d\ge 3$. Let $(f_j)_j$ be an orthonormal system in $L^2_x$.
		For every $\mu\in\ell^2$,
		\begin{equation}\label{eq:C}
			\|\mu_j U(t)f_j\|_{L_t^2\dot B_{\infty,2;x}^{1-d/2}\ell_j^2}
			\lesssim \|\mu\|_{\ell^2}.
		\end{equation}
	\end{lemma}
	\begin{proof}
		The endpoint Strichartz estimate by Keel and Tao \cite{Keel Tao 1998} states that
		\begin{equation}\label{eq:scalar-endpoint}
			\|U(t)f\|_{L_t^2L_x^q}\lesssim \|f\|_{L^2_x},
		\end{equation}
		where $q=\frac{2d}{d-2}$.
		Hence, 
		\begin{align}
			\|\mu_j U(t)f_j\|_{L^2_t \dot{B}^{1-d/2}_{\infty,2;x} \ell_j^2}
			&= \Big\|2^{-(d/2-1)k} \|\mu_j U(t)P_k f_j\|_{L^\infty_x \ell_j^2}\Big\|_{L^2_t\ell_k^2} \\
			&\lesssim \|\mu_j U(t)P_k f_j\|_{ \ell_j^2\ell_k^2 L^2_t L^q_x} \lesssim \|\mu_j \|P_k f_j\|_{\ell_k^2 L^2_x}\|_{\ell_j^2} \lesssim \|\mu\|_{\ell^2}.
		\end{align}
	\end{proof}
	
	\begin{proof}[Proof of Theorem \ref{thm:OC} when $d\ge 3$]
		Let $p\in (1,\infty)$, $s=d/2-1/p$, and $(f_j)_j$ be orthonormal in $\dot{H}^s_x$.
        Let $g_j := |D|^s f_j$.
        Let $s_0 := d/2-1$ and $s_1 := d/2$. Then, $s=s_0/p + s_1/p'$. 
		Then, by Lemma \ref{lem:Besov interpolation}, 
		\begin{align}
			\|U(t)f_j\|_{L^\infty_x \ell_j^2(I)}
            = \||D|^{-s} U(t)g_j\|_{L^\infty_x \ell_j^2(I)}
            \lesssim \|U(t)g_j\|_{\dot B^{-s_0}_{\infty,\infty;x} \ell_j^2(I)}^{1/p} \|U(t)g_j\|_{\dot B^{-s_1}_{\infty,\infty;x} \ell_j^2(I)}^{1-1/p}.
		\end{align}
        Hence, by Lemmas \ref{lem:O} and \ref{lem:C}, 
		\[
		\|U(t)f_j\|_{L^{2p}_t L^\infty_x \ell_j^2(I)}
        \lesssim \|U(t)g_j\|_{L^2_t \dot B^{-s_0}_{\infty,\infty;x} \ell_j^2(I)}^{1/p} \|U(t)g_j\|_{L^\infty_t \dot B^{-s_1}_{\infty,\infty;x} \ell_j^2(I)}^{1-1/p} \lesssim |I|^{\frac1{2p}}.
		\]
		By Lemma \ref{lem:coefficients}, we obtain the desired estimate.
	\end{proof}

	\subsection{When \texorpdfstring{$d=2$}{d=2}}
    We begin with an elementary observation.
	\begin{lemma}\label{lem:local decay}
    Recall that $P_0$ is the Fourier cutoff defined in Section \ref{subsec:LP}. Then,
    \begin{equation}
			\big\|f(x)P_0^2U(t)g(x)\big\|_{\mathfrak{S}^2(L^2_x)} \lesssim \frac{\|f\|_{L^2_x} \|g\|_{L^2_x}}{\langle t \rangle}.
		\end{equation}
	\end{lemma}

    \begin{proof}
    Recall that
    \[
    \|A\|_{\mathfrak S^2(L^2_x)} = \|A(x,y)\|_{L^2_{x,y}},
    \]
    where $A(x,y)$ is the integral kernel of $A$. Hence,
    \begin{equation}
			\big\|f(x)P_0^2U(t)g(x)\big\|_{\mathfrak{S}^2(L^2_x)}
            \le \|f\|_{L^2_x} \|g\|_{L^2_x} \big\|P_0^2U(t)\big\|_{L^\infty_{x,y}}
            \lesssim \frac{\|f\|_{L^2_x} \|g\|_{L^2_x}}{\langle t \rangle}.
	\end{equation}
    \end{proof}

    Next, we give an elementary multilinear estimate.
	\begin{lemma}\label{lem:multilinear 2}
		Let $n\ge 2$ be an integer. For all positive functions $F_j$ and $G_k$, 
		\begin{multline}\label{eq:multilinear 2}
			\int_\R dt_1 \cdots \int_\R dt_n F_1(t_1)\cdots F_n(t_n)G_1(t_1-t_2)\cdots G_n(t_n-t_1) \\
			\le \prod_{j=1}^n \|F_j\|_{L^n_t} \prod_{j=1}^n \|G_j\|_{L^{\frac{n}{n-1}}_t}.
		\end{multline}
	\end{lemma}
	
	\begin{proof}
    For all $\ell,k'\in \{1,\dots,n\}$,
		\begin{multline}
			\sup_{t_\ell \in \R} \int_\R dt_1 \cdots \int_\R dt_{\ell-1} \int_\R dt_{\ell+1} \cdots \int_\R dt_n G_1(t_1-t_2)\cdots G_n(t_n-t_1) \\ 
			\le \|G_{k'}\|_{L^\infty_t} \prod_{k=1,\dots,k'-1,k'+1,\dots,n} \|G_k\|_{L^1_t}.
		\end{multline}
		Hence, by complex interpolation, we obtain
		\begin{multline}\label{eq:multilinear 1}
			\sup_{t_\ell \in \R} \int_\R dt_1 \cdots \int_\R dt_{\ell-1} \int_\R dt_{\ell+1} \cdots \int_\R dt_n G_1(t_1-t_2)\cdots G_n(t_n-t_1)\\
            \le \prod_{k=1}^n \|G_k\|_{L^{\frac{n}{n-1}}_t}
		\end{multline}
        for all $\ell\in \{1,\dots,n\}$.
		By H\"older's inequality,
		\begin{multline}
			\int_\R dt_1 \cdots \int_\R dt_n F_1(t_1)\cdots F_n(t_n)G_1(t_1-t_2)\cdots G_n(t_n-t_1) \\
			\le \|F_1(t_1)\cdots F_n(t_n)\|_{L^1_{t_\ell} L^\infty_{t_1,\dots,t_{\ell-1},t_{\ell+1},\dots,t_n}} \\
			\cdot \|G_1(t_1-t_2)\dots G_n(t_n-t_1)\|_{L^\infty_{t_\ell} L^1_{t_1,\dots,t_{\ell-1},t_{\ell+1},\dots,t_n}}.
		\end{multline}
		Hence, \eqref{eq:multilinear 1} implies
		\begin{multline}
			\int_\R dt_1 \cdots \int_\R dt_n F_1(t_1)\cdots F_n(t_n)G_1(t_1-t_2)\cdots G_n(t_n-t_1) \\
			\le \|F_\ell\|_{L^1_t} \prod_{j=1,\dots,\ell-1,\ell+1,\dots,n} \|F_j\|_{L^\infty_t}
			\prod_k \|G_k\|_{L_t^{\frac{n}{n-1}}}
		\end{multline}
        for every $\ell\in\{1,\ldots,n\}$.
		Finally, by complex interpolation, we obtain \eqref{eq:multilinear 2}.
	\end{proof}
	
	\begin{lemma}\label{lem:localized}
		Let $n\ge2$ be an integer.
		Then, for every $k\in\Z$ and for every $(f_j)_j \subset L^2_x$,
		\begin{equation}\label{eq:localized density}
			\bigg\|\sum_j |U(t)\widetilde P_k f_j|^2\bigg\|_{L_t^{\frac{n}{n-1}} L_x^\infty}
			\lesssim_n 2^{2k/n}\bigg\|\sum_j |f_j\rangle\langle f_j|\bigg\|_{\mathfrak{S}^{\frac{n}{n-1}}(L^2_x)}.
		\end{equation}
		No orthogonality assumption is imposed on this finite family.
	\end{lemma}
	
	\begin{proof}
		By a simple scaling and the triangle inequality, it suffices to prove
        \begin{equation}
			\bigg\|\sum_j |U(t) P_0 f_j|^2\bigg\|_{L_t^{\frac{n}{n-1}}L_x^\infty}
			\lesssim \bigg\|\sum_j |f_j\rangle\langle f_j|\bigg\|_{\mathfrak{S}^{\frac{n}{n-1}}(L^2_x)}.
		\end{equation}
        By the duality principle, it suffices to prove that
		\begin{equation}\label{eq:dual}
			\bigg\|P_0 \int_\R U(t)^* V(t)U(t)dt P_0 \bigg\|_{\mathfrak{S}^n(L^2_x)}
			\lesssim \|V\|_{L_t^n L^1_x}.
		\end{equation}
		Since we can assume $V \ge 0$, 
		\begin{align*}
			&\bigg\|P_0 \int_\R U(t)^* V(t)U(t)dt P_0 \bigg\|_{\mathfrak{S}^n(L^2_x)}^n
			= \Tr \bigg[ \bigg(P_0 \int_\R U(t)^* V(t) U(t) dtP_0\bigg)^n\bigg] \\
			&\qquad = \int_\R dt_1 \cdots\int_\R dt_n \Tr\Big[\sqrt{V(t_1)}U(t_1-t_2)P_0^2\sqrt{V(t_2)} \\
			&\qquad \qquad \qquad \qquad \qquad \qquad \qquad \cdot \cdots \sqrt{V(t_n)} U(t_n-t_1)P_0^2\sqrt{V(t_1)}\Big] \\
			&\qquad \le \int_\R dt_1 \cdots \int_\R dt_n \|\sqrt{V(t_1)}U(t_1-t_2)P_0^2\sqrt{V(t_2)}\|_{\mathfrak{S}^2(L^2_x)}  \\
			&\qquad \qquad \qquad \qquad \qquad \qquad \qquad\cdots \|\sqrt{V(t_n)} U(t_n-t_1)P_0^2 \sqrt{V(t_1)}\|_{\mathfrak{S}^2(L^2_x)} \\
			&\qquad \lesssim \int_\R dt_1 \cdots \int_\R dt_n \frac{\|V(t_1)\|_{L^1_x} \cdots \|V(t_n)\|_{L^1_x}}{\langle t_1-t_2\rangle \cdots \langle t_n - t_1 \rangle}
			\lesssim \|V\|_{L^n_t L^1_x}^n,
		\end{align*}
		where we used Lemma \ref{lem:multilinear 2} to get the last inequality.
	\end{proof}

    \begin{lemma}[Estimate at $(0,1-1/n)$]\label{lem:p_n}
		Let $d=2$, $n\ge 2$, and $p_n := n/(n-1)$.
        Let $(f_j)_j$ be an orthonormal system in $L_x^2$.
		For every $\mu\in\ell^2$,
		\begin{equation}\label{eq:C d=2}
			\|\mu_j U(t)f_j\|_{L_t^{2p_n}\dot B_{\infty,2p_n;x}^{-1+1/p_n}\ell_j^2}
			\lesssim \|\mu\|_{\ell^\infty}^{1/n} \|\mu\|_{\ell^2}^{1-1/n}.
		\end{equation}
	\end{lemma}
    \begin{proof}
    We can assume $\mu_j \ge 0$ for all $j$.
    Clearly, 
		\[
		\bigg\|\sum_j \mu_j^2 |P_k f_j\rangle\langle P_k f_j|\bigg\|_{\mathcal{B}(L^2_x)}
		\le \sup_j \mu_j^2
		\]
		and
		\[
		\bigg\|\sum_j \mu_j^2 |P_k f_j\rangle\langle P_k f_j|\bigg\|_{\mathfrak{S}^1(L^2_x)}
		\le \sum_j\mu_j^2 \|P_k f_j\|_{L^2_x}^2.
		\]
		Therefore, by Lemma \ref{lem:localized},
		\begin{align}
        \|\mu_j U(t)P_k f_j\|_{L^{2p_n}_t L^\infty_x \ell_j^2}
		&\lesssim 2^{k/n} \bigg\|\sum_j \mu_j^2 |P_k f_j \rangle \langle P_k f_j| \bigg\|_{\mathfrak{S}^{p_n}(L^2_x)}^{\frac12} \\
		&\lesssim 2^{k/n} \Big(\sup_j \mu_j\Big)^{1/n}\bigg(\sum_j \mu_j^2 \|P_k f_j\|_{L^2_x}^2 \bigg)^{\frac{1}{2p_n}}.
		\end{align}
		Consequently, the Littlewood--Paley $L^2$ estimate yields
		\[
		\begin{aligned}
			\|\mu_j U(t) f_j\|_{L_t^{2p_n}\dot B_{\infty,2p_n;x}^{-1+1/p_n}\ell_j^2}^{2p_n}
			&\le\sum_{k\in\Z}2^{-\frac{2k}{n-1}}\int_{\R}\|\mu_j U(t)P_k f_j\|_{L_x^\infty\ell_j^2}^{2p_n}dt\\
			&\lesssim \Big(\sup_j \mu_j\Big)^{\frac2{n-1}} \sum_j \mu_j^2 \sum_{k\in\Z}\|P_k f_j\|_{L_x^2}^2 \\
            &\lesssim \Big(\sup_j \mu_j\Big)^{\frac2{n-1}} \sum_j \mu_j^2.
		\end{aligned}
    \]
    \end{proof}

	\begin{proof}[Proof of Theorem \ref{thm:OC} when $d=2$]
		Let $s=1-1/p$, and choose an integer $n\ge 2$ such that $p_n:=n/(n-1)<p$.
		Let $(f_j)_j$ be orthonormal in $\dot{H}_x^s$, and set $g_j=|D|^sf_j$.
		Then $(g_j)_j$ is orthonormal in $L_x^2$.
		Fix a finite set $I\subset\N$.
		By Lemma \ref{lem:p_n}, 
		\[
			\|U(t)g_j\|_{L_t^{2p_n}\dot B_{\infty,2p_n;x}^{-1+1/p_n}\ell_j^2(I)}
			\lesssim |I|^{\frac1{2p_n}}.
		\]
		On the other hand, by Lemma \ref{lem:O},
		\[
		\|U(t)g_j\|_{L_t^\infty\dot B_{\infty,\infty;x}^{-1}\ell_j^2(I)}
		\lesssim 1.
		\]
		Lemma \ref{lem:Besov interpolation} with the above estimates implies
		\[
		\begin{aligned}
			\|U(t)f_j\|_{L_t^{2p}L_x^\infty\ell_j^2(I)}
			&\lesssim \|U(t)g_j\|_{L_t^{2p_n}\dot B_{\infty,\infty;x}^{-1+1/p_n}\ell_j^2(I)}^{p_n/p}
			\|U(t)g_j\|_{L_t^\infty\dot B_{\infty,\infty;x}^{-1}\ell_j^2(I)}^{1-p_n/p}
			\lesssim |I|^{\frac1{2p}}.
		\end{aligned}
		\]
		Therefore, by Lemma \ref{lem:coefficients}, we complete the proof.
	\end{proof}
	
    \section*{Declaration of AI use}
    ChatGPT (GPT-6.1 Sol and GPT-6 Astra) was used to assist in developing the proofs. The authors wrote the entire manuscript from scratch and take full responsibility for its final version.


\begin{thebibliography}{99}
		\bibitem{Bergh Lofstrom 1976 textbook}
		J.~Bergh and J.~L\"ofstr\"om,
		\emph{Interpolation Spaces: An Introduction},
		Springer-Verlag, Berlin-New York, 1976.
		
		\bibitem{Bez et al 2019}
		N.~Bez, Y.~Hong, S.~Lee, S.~Nakamura, and Y.~Sawano,
		\emph{On the Strichartz estimates for orthonormal systems of initial data with regularity},
		Adv. Math. \textbf{354} (2019), 106736.
		
		\bibitem{Bez et al 2024}
		N.~Bez, S.~Kinoshita, and S.~Shiraki,
		\emph{Boundary Strichartz estimates and pointwise convergence
			for orthonormal systems},
		Trans. London Math. Soc. \textbf{11} (2024), e70002.

        \bibitem{Bez et al 2025}
        N. Bez, S. Kinoshita, and S. Shiraki,
        \emph{A note on Strichartz estimates for the wave equation with orthonormal initial data},
        RIMS K\^oky\^uroku Bessatsu B97 (2025), 83--101.
		
		\bibitem{Bez et al 2020}
		N.~Bez, S.~Lee, and S.~Nakamura,
		\emph{Maximal estimates for the Schrödinger equation with orthonormal initial data},
		Selecta Math. (N.S.) \textbf{26} (2020), no. 4, Paper No. 52.

        \bibitem{Carbery 2009}
        A. Carbery,
        \emph{Almost-orthogonality in the Schatten-von Neumann classes},
        J. Operator Theory \textbf{62} (2009), no. 1, 151--158.
        
		\bibitem{Cwikel 1974}
		M.~Cwikel,
		\emph{On $(L^{p_0}(A_0),L^{p_1}(A_1))_{\theta,q}$},
		Proc. Amer. Math. Soc. \textbf{44} (1974), 286--292.
   		
		\bibitem{Feng et al 2026-}
		G. Feng, J. B. Lee, and S. Lee,
		\emph{Orthonormal Strichartz estimates over fractal time},
		arXiv preprint (2026), arXiv:2610.18042v1.
		
		\bibitem{Frank et al 2014}
		R. L.~Frank, M. Lewin, E. H.~Lieb, and R.~Seiringer,
		\emph{Strichartz inequality for orthonormal functions},
		J. Eur. Math. Soc. \textbf{16} (2014), no. 7, 1507--1526.
		
		\bibitem{Frank Sabin 2016}
		R. L. Frank and J. Sabin,
		\emph{The Stein-Tomas inequality in trace ideals},
		S\'eminaire Laurent Schwartz—\'Equations aux d\'eriv\'ees partielles et applications. Ann\'ee 2015–2016, Exp. No. XV, 12 pp.
		
		\bibitem{Frank Sabin 2017}
		R. L. Frank and J. Sabin,
		\emph{Restriction theorems for orthonormal functions, Strichartz inequalities, and uniform Sobolev estimates},
		Amer. J. Math. \textbf{139} (2017), 1649--1691.
		
		\bibitem{Ginibre Velo 1985}
		J.~Ginibre and G.~Velo,
		\emph{The global Cauchy problem for the non linear Schrödinger equation revisited},
		Annales de l’Institut Henri Poincar\'e, Analyse non lin\'eaire \textbf{2} (1985), 309--327.
				
		\bibitem{Ginibre Velo 1992}
		J.~Ginibre and G.~Velo,
		\emph{Smoothing properties and retarded estimates for some dispersive evolution equations},
		Comm. Math. Phys. \textbf{144} (1992), no.~1, 163--188.
				
		\bibitem{Gong 2026-}
		Y. Gong,
		\emph{Critical Strichartz estimates for orthonormal systems},
		arXiv preprint (2026), arXiv:2608.22805v2.
		
		\bibitem{Gong 2026- 2}
		Y. Gong,
		\emph{Critical orthonormal Strichartz estimates for Lebesgue spaces},
		arXiv preprint (2026), arXiv:2611.03527v2.
		
		\bibitem{Guo et al 2018}
		Z. Guo, J. Li, K. Nakanishi, and L. Yan,
		\emph{On the boundary Strichartz estimates for wave and Schrödinger equations},
		J. Differential Equations \textbf{265} (2018), no. 11, 5656--5675.
				
		\bibitem{Hadama Toyoshima 2026+}
		S.~Hadama and S.~Toyoshima,
		\emph{Endpoint orthonormal Strichartz estimates},
		arXiv preprint (2026), arXiv:2608.21823v1.
		
		\bibitem{Hernandez Soria 1991}
		E.~Hernández and J.~Soria,
		\emph{Spaces of Lorentz type and complex interpolation},
		Arkiv f\"or Matematik \textbf{29} (1991), 203--220.
		
		\bibitem{Keel Tao 1998}
		M.~Keel and T.~Tao,
		\emph{Endpoint Strichartz estimates},
		Amer. J. Math. \textbf{120} (1998), 955--980.
		
		\bibitem{Lieb 1983}
		E. H.~Lieb,
		\emph{An $L^p$ bound for the Riesz and Bessel potentials of orthonormal functions},
		J. Funct. Anal. \textbf{51} (1983), no. 2, 159--165.

        \bibitem{Lindemulder Lorist 2022}
        N. Lindemulder and E. Lorist,
        \emph{Stein interpolation for the real interpolation method},
        Banach J. Math. Anal. \textbf{16} (2022), Paper No.~7.
        
		\bibitem{Lions Peetre 1964}
		J.-L.~Lions and J. Peetre,
		\emph{Sur une classe d'espaces d'interpolation},
		Inst. Hautes \'Etudes Sci. Publ. Math. No. \textbf{19} (1964), 5--68.
		
		\bibitem{Montgomery-Smith 1998}
		S. J. Montgomery-Smith, 
		\emph{Time decay for the bounded mean oscillation of solutions of the Schr\"odinger and wave equations}, 
		Duke Math. J. \textbf{91} (1998), no. 2, 393--408.

        \bibitem{Stein 1993 textbook}
        E. M. Stein,
        \emph{Harmonic analysis: real-variable methods, orthogonality, and oscillatory integrals},
        Princeton University Press, Princeton, NJ, 1993. xiv+695 pp.
        
		\bibitem{Strichartz 1977}
		R. S. Strichartz,
		\emph{Restrictions of Fourier transforms to quadratic surfaces and decay of solutions of wave equations},
		Duke Math. J. \textbf{44} (1977), 705--714.
		
		\bibitem{Tao 2006}
		T. Tao,
		\emph{A counterexample to an endpoint bilinear Strichartz inequality},
		Electron. J. Differential Equations 2006, No. 151, 6 pp.
	\end{thebibliography}
\end{document}